\documentclass[a4paper]{scrartcl} %
\usepackage{etex}

\PassOptionsToPackage{x11names}{xcolor}

\usepackage[utf8]{inputenc}
\usepackage[T1]{fontenc}
\usepackage[english]{babel}
\usepackage[neverdecrease]{paralist}
\usepackage{amssymb}
\usepackage{stmaryrd}
\usepackage[ntheorem]{empheq} %
\usepackage[thmmarks,amsmath]{ntheorem}
\usepackage[pdfusetitle,    colorlinks,
    linkcolor={red!70!black},
    citecolor={green!45!black},
    urlcolor={blue}]{hyperref}  %
\usepackage{jensntheorem}
\usepackage{lmodern}
\usepackage{microtype}
\usepackage{verbatim}
\usepackage{mathrsfs}
\usepackage{graphicx}
\usepackage{dsfont}
\usepackage{tikz-cd}
\usepackage{mathdots}
\usepackage[normalem]{ulem}
\usepackage{lscape}  
\usepackage{geometry}

\numberwithin{equation}{section}

\usepackage{csquotes}
\usepackage[style=alphabetic,maxnames=10,maxalphanames=5,backend=biber]{biblatex}
\bibliography{prom}

\usetikzlibrary{arrows,positioning,decorations.markings}

\usepackage{jensstuff}

\title{EKOR strata on Shimura varieties with parahoric reduction}
\author{Jens Hesse
\thanks{Technische Universit\"at Darmstadt, \texttt{math@jenshesse.eu}}}
\date{\today}

\begin{document}

\maketitle
	
\begin{abstract}
  We investigate the geometry of the special fiber of the integral model of a Shimura variety with parahoric level at a given prime place.

To be more precise, we deal with the EKOR stratification which interpolates between the Ekedahl-Oort and Kottwitz-Rapoport stratifications. In the Siegel case we give a geometric description by suitably generalizing the theory of $G$-zips of Moonen, Wedhorn, Pink and Ziegler to our context.
\end{abstract}

\tableofcontents

\addsec{Introduction}

Shimura varieties are objects of arithmetic geometry (namely varieties over number fields) that naturally arise in the search for generalized, non-abelian reciprocity laws (i.e., in the Langlands program) and as moduli spaces of abelian varieties (with certain extra structures on them). One way of approaching these objects is to try to understand their mod-$p$ reduction (which has to be carefully defined first). Insofar as a moduli interpretation in the above sense exists and continues to exist likewise for the mod-$p$ reduction\footnote{There need not be a \emph{literal} moduli interpretation, but in any event the stratifications in question derive from a close connection to moduli problems.}, it allows us to stratify the moduli space according to  several invariants of the abelian varieties parametrized, e.g., the isomorphism classes of their $p$-torsion. (An important observation is that these stratifications genuinely live in the characteristic $p$ world, making use of Frobenius endomorphisms and so on.) This, very roughly, is the general theme everything in this article revolves around.

More precisely, we will be dealing with Shimura varieties of Hodge type and parahoric level structure, at some fixed prime $v \mid p$ of the number field over which the Shimura variety is defined. Under some reasonably mild assumptions, cf.~\ref{std-assum}, Kisin and Pappas \cite{kisin-pappas} constructed a canonical integral model for such a Shimura variety. We try to understand some aspects of the geometry of the special fiber  of said integral model, namely the EKOR strata (an interpolation between the Ekedahl-Oort strata, which in the case of hyperspecial level are roughly the patches where the isomorphism class of the $p$-torsion associated with the abelian variety is constant, and the Kottwitz-Rapoport strata, which roughly are the patches where the Hodge filtration looks constant) and defining them in a geometrical way.

Let us now go into more detail.

On the integral model $\sS_K$ ($K$ parahoric level) we have a ``universal'' abelian scheme (the quotation marks indicating that it is not really universal for some moduli problem on $\sS_K$, but it comes from a universal abelian scheme via pullback) and we have various kinds of Hodge tensors. We also have a ``universal'' isogeny chain of abelian schemes tightly connected to the ``universal'' abelian scheme. 

The overarching goal (and what we meant above by ``defining the EKOR strata in a geometrical way'') is to construct a ``nice'' algebraic stack $\ocG_K\EKORZip$ and a ``nice'' morphism $\osS_K\to\ocG_K\EKORZip$ from the mod-$p$ reduction of the Shimura variety to it, such that the fibers are the EKOR strata. Shen, Yu and Zhang \cite{SYZnew} solved this problem on individual Kottwitz-Rapoport strata and globally after perfection, but not in the form stated here (i.e., globally without passing to perfections). In the Siegel case we propose a solution which specializes to that of Shen, Yu and Zhang on Kottwitz-Rapoport strata, and should not be difficult to generalize to many (P)EL cases. We show that $\osS_K\to\ocG_K\EKORZip$ is surjective. However, we have to leave the question of whether $\osS_K\to\ocG_K\EKORZip$ is smooth (which would be part of ``nice'') an open conjecture.

For hyperspecial level, the EKOR stratification agrees with the Ekedahl-Oort stratification, and the goal just set out is achieved by the stack of $\ocG_K$-zips, first defined in special cases by Moonen and Wedhorn in \cite{mw} and then generally by Pink, Wedhorn and Ziegler in \cite{PWZ-AZD,PWZ}; the relation to Shimura varieties being established in increasing generality in \cite{mw}, by Viehmann and Wedhorn in \cite{vieh-wed}, and finally by Zhang in \cite{zhang}.

One way of looking at the transition from hyperspecial level to general parahoric level (at the very least in nice enough (P)EL cases) is from the point of view of moduli problems of abelian varieties with extra structure, where in the hyperspecial case we are really dealing just with that and in the general case we are dealing with isogeny chains of abelian varieties with extra structure, indexed by lattice chains coming from the Bruhat-Tits building of the reductive $p$-adic Lie group in question.
The basic idea in generalizing zips from the hyperspecial to the general parahoric case then is that one should be dealing with chains of zips in the old sense.

The zip of an abelian variety encodes the following information: the Hodge filtration, the conjugate filtration, and the Cartier isomorphism relating the two. In the general case, every abelian variety in the isogeny chain has a Hodge filtration, a conjugate filtration and a Cartier isomorphism. Problems now arise because we are dealing with $p$-primary isogenies on $p$-torsion points, implying that the transition morphisms in these chains have non-vanishing kernels.  This introduces additional difficulty compared to the hyperspecial case; there is a naive way of defining a zip stack, but eventually we need to consider a certain admissible locus in it, which so far suffers from the absence of a nice moduli description.  Passing to perfections however simplifies things and allows us to prove that the admissible locus is closed.  From here we arrive at the stack that we are really interested in by dividing out a certain group action involving the unipotent radical of the special fiber of the parahoric group scheme.  A careful inspection shows that on Kottwitz-Rapoport strata we arrive at the same result  as in \cite{SYZnew}.

To sum up the results,
\begin{TheoremA}
  In the Siegel case, there is an algebraic stack $\ocG_K\EKORZip$ and a        surjective         morphism $\osS_K \to \ocG_K\EKORZip$, whose fibers are the EKOR strata and such that on Kottwitz-Rapoport strata, one gets the stack and map constructed in \cite{SYZnew}.
\end{TheoremA}

For $\GSp(4)$ we do some calculations to illustrate the theory; section~\ref{sec:example-gsp4}.

\subsection*{Acknowledgements}

This article essentially is an extract of my doctoral thesis \cite{diss} (another extract\footnote{In particular, there is a large overlap between the ``Background'' sections of the two articles.}, dealing with the foliation into central leaves, is \cite{article-leaves}). I thank Torsten Wedhorn for suggesting the topic of the dissertation, his support and patiently answering my questions.
Moreover I thank Eva Viehmann and Paul Hamacher for their hospitality and helpful discussions during a month-long stay in Munich at the TU M\"unchen.
I am also grateful to Timo Richarz and Timo Henkel
for numerous helpful discussions.

This research was supported by the Deutsche Forschungsgemeinschaft (DFG), project number WE~2380/5.

\section{Background}
\label{cha:gener-prel}

\subsection{Shimura data of Hodge type}
\label{sec:shimura-data-hodge}

This article deals with aspects of the geometry of Shimura varieties (of Hodge type), which are the (systems of) varieties associated with Shimura data (of Hodge type).

\begin{Definition}\label{def-hodget}
  A \defn{Shimura datum of Hodge type}\index{Shimura datum of Hodge type} is a pair $(G,X)$, where $G$ is a reductive algebraic group over $\Q$ and $X\subseteq\Hom_{\R\text{-grp}}(\SSS,G_\R)$ is a $G(\R)$-conjugacy class ($\SSS:=\Res_{\C/\R}\G_{m,\C}$ being the Deligne torus\index{Deligne torus}) subject to the following conditions:
\begin{enumerate}[(1)]
\item For $h\in X$, the induced Hodge structure $\SSS\xrightarrow{h} G_\R \xrightarrow{\mathrm{Ad}}\GL(\Lie(G_\R))$ satisfies $\Lie(G_\C)=\Lie(G_\C)^{-1,1}\oplus\Lie(G_\C)^{0,0}\oplus\Lie(G_\C)^{1,-1}$.\label{sv1}
\item $\interior(h(i))\colon G^\mathrm{ad}_\R\to G^\mathrm{ad}_\R$ is a Cartan involution, i.e., $\{ g\in G^\mathrm{ad}(\C)\suchthat gh(i) = h(i)\overline{g}\}$ is compact. Another way of phrasing this condition: Every finite-dimensional real representation $V$ of $G^\mathrm{ad}_\R$ carries a $G^\mathrm{ad}_\R$-invariant bilinear form $\varphi$ such that $(u,v)\mapsto \varphi(u,h(i)v)$ is symmetric and positive definite. It is enough to show that this holds for one \emph{faithful} finite-dimensional real representation $V$.
\item $G^\mathrm{ad}$ \emph{cannot} be non-trivially written as $G^\mathrm{ad}\cong H\times I$ over $\Q$ with $\SSS\to G_\R \xrightarrow{\mathrm{proj}} H_\R$ trivial.
\item There exists an embedding $(G,X)\hookrightarrow (\GSp(V),S^\pm)$, where $(\GSp(V),S^\pm)$ is the Shimura datum associated with a finite-dimensional symplectic $\Q$-vector space $V$ (see below). That is, we have an embedding $G\hookrightarrow \GSp(V)$ of $\Q$-group schemes such that the induced map $\Hom_{\R\text{-grp}}(\SSS,G_\R)\hookrightarrow\Hom_{\R\text{-grp}}(\SSS,\GSp(V_\R))$ restricts to a map $X\hookrightarrow S^\pm$. \label{item:def-hodget-hodge-emb}
\end{enumerate}
\end{Definition}

\begin{Example}\label{sympl-ex}
  Let $W$ be a finite-dimensional $\R$-vector space.

  $\R$-group homomorphisms $\SSS\to \GL(W)$ then correspond to Hodge decompositions of $W$, i.e., to decompositions $W_\C=\oplus_{(p,q)\in\Z^2}W_\C^{p,q}$, such that $W_\C^{p,q}$ is the complex conjugate of $W_\C^{q,p}$ for all $(p,q)\in\Z^2$. Under this correspondence, $h\colon \SSS\to \GL(W)$ corresponds to the Hodge decomposition $W_\C^{p,q}=\{ w\in W_\C \suchthat \forall z\in\SSS(\R)=\C^\times\colon h(z)w=z^{-p}\bar{z}^{-q}w \}$. Hodge decompositions of $W$ of type $(-1,0)+(0,-1)$ correspond to complex structures on $W$: If $h\colon\SSS\to\GL(W)$ yields such a Hodge decomposition, then $h(i)$ gives an $\R$-endomorphism $J$ of $W$ with $J\circ J=-\id_W$.

  Let $V=(V,\psi)$ be a finite-dimensional symplectic $\Q$-vector space. We say that a complex structure $J$ on $V_\R$ is positive (resp. negative) if $\psi_J:=\psi_\R(\_,J\_)$ is a positive definite (resp. negative definite) symmetric bilinear form on $V_\R$. Define $S^+$ (resp. $S^-$) to be the set of positive (resp. negative) complex structures on $(V_\R,\psi_\R)$ and $S^\pm:=S^+\sqcup S^-$.

   We can make this more concrete: A symplectic basis of $(V_\R,\psi_\R)$ is a basis $e_1,\dotsc,\allowbreak e_g, e_{-g}, \dotsc, \allowbreak e_{-1}$, such that $\psi_\R$ is of the form
  \begin{equation*}
    \begin{pmatrix}
      & \tilde{I}_g \\
      -\tilde{I}_g &
    \end{pmatrix}
  \end{equation*}
  with respect to this basis, where $\tilde{I}_g=
  \begin{pmatrix}
    & & 1\\
    & \iddots & \\
    1 & &
  \end{pmatrix}$ is the antidiagonal identity matrix.\footnote{Occasionally (in particular when doing concrete matrix calculations), it is more convenient to number the basis vectors $1,\dotsc,g,-1,\dotsc,-g$ instead of $1,\dotsc,g,-g,\dotsc,-1$. Then the standard symplectic form is given by $\begin{smallpmatrix}
      & I_g \\
      -I_g &
    \end{smallpmatrix}$, $I_g$ being the $g\times g$ identity matrix.}

  Let $J$ be the endomorphism of $V_\R$ of the form
  \begin{equation*}
    \begin{pmatrix}
      & -\tilde{I}_g \\
      \tilde{I}_g &
    \end{pmatrix}
  \end{equation*}
  with respect to this basis. Then $J\in S^+$ and what we have described is a surjective map
  \begin{equation*}
    \{\text{symplectic bases of } (V_\R,\psi_\R)\} \twoheadrightarrow S^+.
  \end{equation*}

  In particular we see that $\Sp(V_\R,\psi_\R):=\{f\in\GL(V_\R) \suchthat \psi_\R(f(\_),f(\_))=\psi_\R\}$  (by virtue of acting simply transitively on the symplectic bases) acts transitively on $S^+\cong \Sp(V_\R,\psi_\R)/\SpO(V_\R,\psi_\R,J)$ (where we define $\SpO(V_\R,\psi_\R,J):=\Sp(V_\R,\psi_\R)\cap O(V_\R,\psi_J)=U((V_\R,J),\psi_J)$ for a fixed choice of $J\in S^+$) and therefore the general symplectic group $\GSp(V_\R,\psi_\R):=\{f\in\GL(V_\R) \suchthat \psi_\R(f(\_),f(\_))=c\cdot\psi_\R\text{ for some } c\in\R^\times\}$ acts transitively on $S^\pm$ (note that the element of the form $e_{\pm i}\mapsto e_{\mp i}$ of $\GSp(V_\R,\psi_\R)$ for any given choice of symplectic basis $\left(e_i\right)_i$ permutes $S^+$ and $S^-$).
\end{Example}

\begin{Definition}
  Condition~\eqref{sv1} of Definition~\ref{def-hodget} implies that the action of $\G_{m,\R}$ (embedded in $\SSS$ in the natural way) on $\Lie(G_\R)$ is trivial, so that $h$ induces a homomorphism ${w\colon \G_{m,\R}\to \Cent(G_\R)}$. This homomorphism is independent of the choice of $h\in X$ and is called the \defn{weight homomorphism} of $(G,X)$.

  Moreover, we denote by $\{\mu\}$ the the $G(\C)$-conjugacy class of the cocharacter $\mu_h:=h\circ(\id_{\G_{m,\C}},1)\colon \G_{m,\C}\to \G_{m,\C}^2\cong\SSS_\C\to G_\C$, where $h$ is as above. Obviously, the conjugacy class $\{\mu\}$ is independent of the particular choice of $h\in X$.
\end{Definition}

\begin{Remark}\label{conj-class-cochar}
  Let $L/\Q$ be a field extension such that $G_L$ contains a split maximal torus $T$.  Let $W:=\Norm_{G(L)}(T)/T$ be the Weyl group.  Then the natural map
  \begin{equation*}
    W\backslash \Hom_{L\text{-grp}}(\G_{m,L},T) \to G(L)\backslash \Hom_{L\text{-grp}}(\G_{m,L},G_L)
  \end{equation*}
  is bijective.

  Since the left hand side remains unchanged if we go from $L=\bar\Q$ (where as usual $\bar\Q$ denotes an algebraic closure of $\Q$) to $L=\C$, we see that $\{\mu\}$ contains a cocharacter defined over $\bar\Q$ and that we may then also consider $\{\mu\}$ as a $G(\bar\Q)$-conjugacy class.
\end{Remark}

\begin{Definition}
  The \defn{reflex field} $\bfE=\bfE(G,X)$ of $(G,X)$ is the field of definition of $\{\mu\}$, i.e., the fixed field in $\bar\Q$ of $\{\gamma\in\Gal(\bar\Q/\Q) \suchthat \gamma(\{\mu\})=\{\mu\}\}$.
\end{Definition}

\begin{Example}
  The reflex field of the Shimura datum $(\GSp_{2g,\Q},S^\pm)$ of Example~\ref{sympl-ex} is $\Q$. To wit, one of the cocharacters in the conjugacy class $\{\mu\}$ is
  \begin{equation*}
    \mu(z) =
    \begin{smallpmatrix}
      z & & & & & \\
    & \ddots & & & &\\
    & & z & & &\\
    & & & 1 & &  \\
    & & & & \ddots & \\
    & & & & & 1
    \end{smallpmatrix}.
  \end{equation*}
\end{Example}

\begin{Notation}
We denote the ring of (rational) adeles by $\A:=\A_\Q$, the subring of finite adeles by $\A_f:=\A_{\Q,f}$ and the subring of finite adeles away from some fixed prime $p$ by $\A_f^p$.
\end{Notation}

\begin{DefinitionRemark}\label{shvar}
  Let $K\subseteq G(\A_f)$ be a compact open subgroup. The \defn{Shimura variety of level $K$ associated with $(G,X)$} is the double coset space
  \begin{equation*}
    \Sh_K(G,X):=G(\Q)\backslash (X\times (G(\A_f)/K)).
  \end{equation*}

  A priori, this is just a set, but if $K$ is sufficiently small (i.e., ``neat'' in the sense of \cite{borelarith,pink}), $\Sh_K(G,X)$ can be canonically written as a finite disjoint union of hermitian symmetric domains.\footnote{If $K$ fails to be sufficiently small, one might very reasonably argue that our definition of the Shimura variety of level $K$ really is the definition of the \emph{coarse} Shimura variety and that one should be working with stacks instead.  Since we will only be interested in sufficiently small level, this is inconsequential for us.}  In particular, this gives $\Sh_K(G,X)$ the structure of a complex manifold.

  In fact, by the theorem of Baily-Borel, this complex manifold attains the structure of a quasi-projective complex variety in a canonical way.

  By work of Deligne, Milne and Borovoi, this variety is defined already (and again in a canonical way) over the reflex field $\bfE$. So in particular, it is defined over a number field independent of $K$.  This is important when varying $K$ and it is the reason why we consider the whole Shimura variety instead of its connected components over $\C$ on their own.  It is possible for the Shimura variety to have multiple connected components over $\C$ while being connected over $\bfE$.

  More detailed explanations may be found in \cite{milne-isv}.
\end{DefinitionRemark}

\subsection{Bruhat-Tits buildings}
\label{sec:bruh-tits-build}

Let $K$ be a complete discrete valuation field with ring of integers $\mathcal{O}$, uniformizer $\varpi$ and perfect residue field $\kappa:=\mathcal{O}/\varpi$.

\begin{Notation}
  For a (connected) reductive group $G$ over $K$, we denote by $\mathcal{B}(G,K)$ the extended (or enlarged) and by $\mathcal{B}^\mathrm{red}(G,K)$ the reduced (i.e., non-extended) Bruhat-Tits building of $G$ over $K$ \cite{bt-ii}. Moreover, $\mathcal{B}^\mathrm{abstract}(G,K)$ denotes the underlying abstract simplicial complex.
\end{Notation}

\begin{Remark}\label{gl-building}
  Let $V$ be a finite-dimensional $K$-vector space.
  
  As described in \cite[\nopp 1.1.9]{kisin-pappas} (originally in \cite{zbMATH03900941}), the points of $\mathcal{B}(\GL(V),K)$ correspond to graded periodic lattice chains $(\mathcal{L},c)$, i.e.,
\begin{itemize}
\item $\emptyset\ne\mathcal{L}$ is a totally ordered set of full $\mathcal{O}$-lattices in $V$ stable under scalar multiplication (i.e., $\Lambda\in\mathcal{L} \iff \varpi\Lambda\in\mathcal{L}$),
\item $c\colon \mathcal{L}\to\R$ is a strictly decreasing function such that $c(\varpi^n\Lambda)=c(\Lambda)+n$.
\end{itemize}
\end{Remark}

\begin{Remark}\label{gitterketten-nummerieren}
  Fix such an $\mathcal{L}$ and let $\Lambda^0\in\mathcal{L}$. Then every homothety class of lattices has a unique representative $\Lambda$ such that $\Lambda\subseteq\Lambda^0$ and $\Lambda\not\subseteq \varpi\Lambda^0$. Consider such representatives $\Lambda^i$ for all of the distinct homothety classes of lattices that make up $\mathcal{L}$. Because $\mathcal{L}$ is totally ordered and $\Lambda^i\not\subseteq \varpi\Lambda^0$, it follows that $\Lambda^i\supseteq \varpi\Lambda^0$ for all $i$ and that $\left\{\Lambda^i/\varpi\Lambda^0\right\}_i$ is a flag of non-trivial linear subspaces of $\Lambda^0/\varpi\Lambda^0\cong\kappa^{n}$, where $n:=\dim V$. Consequently, the number $r$ of homothety classes is in $\{1,\dotsc,n\}$; it is called the \defn{period length} (or \defn{rank}) of $\mathcal{L}$. Numbering the $\Lambda^i$ in descending order we hence obtain $r$ lattices $\Lambda^0,\Lambda^1,\dotsc,\Lambda^{r-1}$ such that
  \begin{equation}\label{eq:numbered-lattice-chain}
    \Lambda^0 \supsetneqq \Lambda^1 \supsetneqq \dotsb\supsetneqq \Lambda^{r-1} \supsetneqq \varpi\Lambda^0
  \end{equation}
  and $\mathcal{L}$ is given by the the strictly descending sequence of lattices
  \begin{equation*}
    \Lambda^{qr+i}=\varpi^q\Lambda^i,\quad q\in\Z, \; 0\leq i < r.
  \end{equation*}
\end{Remark}

\begin{Remark}\label{gsp-building}
  Let $V$ be a finite-dimensional symplectic $K$-vector space.
  
  $\mathcal{B}(\GSp(V),K)$ embeds into the subset of $\mathcal{B}(\GL(V),K)$ consisting of those $(\mathcal{L},c)$ such that
  $\Lambda\in\mathcal{L}\implies \Lambda^\vee\in \mathcal{L}$.

  Passing to the underlying abstract simplicial complexes means forgetting about the grading $c$ and
  \begin{equation*}
    \mathcal{B}^\mathrm{abstract}(\GSp(V),K) = \{\mathcal{L}\in\mathcal{B}^\mathrm{abstract}(\GL(V),K) \suchthat \Lambda\in\mathcal{L}\implies \Lambda^\vee\in \mathcal{L}\}.
  \end{equation*}

  If $\mathcal{L}\in\mathcal{B}^\mathrm{abstract}(\GSp(V),K)$ and $\{\Lambda^i\}_i$ is as in Remark~\ref{gitterketten-nummerieren}, then there is an involution $t\colon \Z\to\Z$ with $\left(\Lambda^i\right)^\vee=\Lambda^{t(i)}$, $t(i+qr)=t(i)-qr$, and $i< j\implies t(i)> t(j)$. So $-a:=t(0)> t(1)> \dotsb > t(r)=-a-r$, which implies $t(i)=-i-a$. Thus $i_0-t(i_0)=2i_0+a\in\{0,1\}$ for some unique $i_0\in\Z$. Hence, upon renumbering the $\Lambda^i$, we may assume that $a\in\{0,1\}$.

We therefore have
\begin{align*}
  \varpi\Lambda^0\subsetneqq\Lambda^{r-1}\subsetneqq \Lambda^{r-2}\subsetneqq \dotsb \subsetneqq \Lambda^0 \subseteq \left(\Lambda^{0}\right)^\vee=\Lambda^{-a} \subsetneqq \left(\Lambda^{1}\right)^\vee=\Lambda^{-1-a}  \\
  \subsetneqq \dotsb\subsetneqq \left(\Lambda^{r-1}\right)^\vee=\Lambda^{-r+1-a} \subseteq \Lambda^{-r}=\varpi^{-1}\Lambda^0.
\end{align*}
\end{Remark}

\begin{Example}
  See also section~\ref{sec:example-gsp4} for some elaborations on the building of $\GSp_{4}(\Q_p)$.
\end{Example}

\subsection{Bruhat-Tits group schemes}
\label{sec:bruhat-tits-group}

\begin{Notation}
  Let $E$ be a finite field extension of $\Q_p$.
  
  Denote by $\bE$ the completion of the maximal unramified extension of $E$ (hence $\bE=E\cdot\bQ_p$).
\end{Notation}

\begin{Remark}\label{ram-witt}
  If $E/\Q_p$ is unramified, then $\bOE=W(\bar\F_p)$, $\bar\F_p$ denoting an algebraic closure of $\F_p$ and $W\colon\mathrm{Ring}\to\mathrm{Ring}$ being the ($p$-adic) Witt vectors functor. This generalizes to the ramified case using \emph{ramified Witt vectors} instead, see e.g. \cite[Chap.~IV, (18.6.13)]{haze-fg} or \cite[Chapter~1]{ahsendorf}.
\end{Remark}

Let $(G,X)$ be a Shimura datum of Hodge type, let $(G,X)\hookrightarrow(\GSp(V),S^\pm)$ be an embedding as in Definition~\ref{def-hodget}\,\eqref{item:def-hodget-hodge-emb}, and let $x\in\mathcal{B}(G,\Q_p)$ be a point in the Bruhat-Tits building of $G$ over $\Q_p$.

We consider the associated Bruhat-Tits scheme ${\cal G}_x$, i.e., the affine smooth model of $G_{\Q_p}$ over $\Z_p$ such that ${\cal G}_x(\bZ_p)\subseteq G(\bQ_p)$ is the stabilizer of the facet of $x$ in ${\cal B}(G,\bQ_p)\overset{\text{\cite[Prop.~2.1.3]{landvogt}}}=\mathcal{B}(G,\uQ_p)$.
Let $K_p:={\cal G}_x(\Z_p)\subseteq G(\Q_p)$ and let $K^p\subseteq G(\A_f^p)$ be a sufficiently small open compact subgroup. Define $K:=K_pK^p\subseteq G(\A_f)$.

\begin{Assumptions}\label{std-assum}
  From now on, we will always make the following assumptions:
  \begin{itemize}
  \item $\mathcal{G}_x=\mathcal{G}_x^\circ$ is connected.
  \item $G$ splits over a tamely ramified extension of $\Q_p$.
  \item $p\nmid \#\pi_1(G^\mathrm{der})$.
  \end{itemize}
\end{Assumptions}

\begin{Notation}
  In order not to make notation overly cumbersome, we usually denote the base change $G_{\Q_p}$ of $G$ to $\Q_p$ by $G$ again.  (Later, we will almost exclusively be dealing with $G_{\Q_p}$.)
\end{Notation}

\subsection{Siegel integral models}
\label{sec:siegel-integr-model}

With notation as above let
\begin{align*}
  N_p &:= \Stab_{\GSp(V)(\Q_p)}(\mathcal{L}) \quad\text{(as before)}, \\
  J_p&:= \Stab_{\GL(V^\S)(\Q_p)}(\Lambda^\S)\cap\GSp(V^\S)(\Q_p).
\end{align*}
Let $N^p\subseteq\GSp(V)(\A_f^p)$ and $J^p\subseteq\GSp(V^\S)(\A_f^p)$ be sufficiently small open compact subgroups, and $N:=N_pN^p$, $J:=J_pJ^p$.

In this subsection, we are going to describe integral models of $\Sh_{N}(\GSp(V),S^\pm)$ and of $\Sh_{J}(\GSp(V^\S),S^{\S,\pm})$ over $\Z_{(p)}$ and relate the two.

\begin{Remark}\label{rz-moduli}
  By \cite[Definition~6.9]{rz}, the integral model $\sS_{N}(\GSp(V),S^\pm)$ is given by the moduli problem $(\Z_{(p)}\text{-scheme})\ni S\mapsto \left\{(A,\bar\lambda,\eta^p)\right\}/\littlecong$, where:
  \begin{enumerate}[(a)]
  \item $A=\left(A_\Lambda\right)_{\Lambda\in\mathcal{L}}$ is an $\mathcal{L}$-set of abelian schemes, i.e.,
    \begin{itemize}
    \item for every $\Lambda\in\mathcal{L}$, an abelian $S$-scheme up to $\Z_{(p)}$-isogeny $A_\Lambda$ (i.e., $A_\Lambda$ is an object of the category $(\text{abelian } S\text{-schemes})\otimes\Z_{(p)}$, where the category $\mathcal{A}\otimes R$ for $\mathcal{A}$ an preadditive category and $R$ a ring has the same objects as $\mathcal{A}$ and $\Hom_{\mathcal{A}\otimes R}(X,Y)=\Hom(X,Y)\otimes_\Z R$ for all objects $X,Y$),
    \item for every inclusion $\Lambda_1\subseteq \Lambda_2$ a $\Z_{(p)}$-isogeny $\rho_{\Lambda_2,\Lambda_1}\colon  A_{\Lambda_1}\to A_{\Lambda_2}$,
    \item $\rho_{\Lambda_3,\Lambda_1}=\rho_{\Lambda_3,\Lambda_2}\circ\rho_{\Lambda_2,\Lambda_1}$ if $\Lambda_1\subseteq \Lambda_2 \subseteq \Lambda_3$ in $\mathcal{L}$,
    \item the height of $\rho_{\Lambda_2,\Lambda_1}$ is $\log_p|\Lambda_2/\Lambda_1|$. Here $\rho_{\Lambda_2,\Lambda_1}$ gives rise to a well-defined homomorphism of $p$-divisible groups, and what we mean is that the kernel of this homomorphism (which is a finite locally free commutative group scheme, which we also refer to simply as the kernel of $\rho_{\Lambda_2,\Lambda_1}$) is to have order $|\Lambda_2/\Lambda_1|$.
    \item For every $\Lambda\in \mathcal{L}$, there is an isomorphism (called \defn{periodicity isomorphism}) $\theta_\Lambda\colon A_\Lambda\to A_{p\Lambda}$ such that $\rho_{\Lambda,p\Lambda}\circ \theta_\Lambda = [p]\colon A_\Lambda\to A_\Lambda$ is the multiplication-by-$p$ isogeny.
    \end{itemize}
  \item $\bar\lambda\colon A\to\tilde{A}$ is a $\Q$-homogeneous principal polarization, i.e., a $\underline{\Q^\times}$-orbit of a principal polarization $\lambda\colon A\to \tilde{A}$. Here $\tilde{A}$ is the $\mathcal{L}$-set of abelian schemes over $S$ up to prime-to-$p$ isogeny given by $\tilde{A}_\Lambda:=(A_{\Lambda^\vee})^\vee$. And being a polarization $\lambda$ means being a quasi-isogeny of $\mathcal{L}$-sets $\lambda\colon A\to\tilde{A}$ such that
    \begin{equation*}
      A_\Lambda \xrightarrow{\lambda_\Lambda}\tilde{A}_\Lambda=(A_{\Lambda^\vee})^\vee\xrightarrow{\varrho_{\Lambda^\vee,\Lambda}^\vee}(A_\Lambda)^\vee
    \end{equation*}
    is a polarization of $A_\Lambda$ for all $\Lambda$. If $\lambda_\Lambda$ can be chosen to be an isomorphism up to prime-to-$p$ isogeny for all $\Lambda$, then we speak of a principal polarization. In that case, when referring to $\lambda_\Lambda$, we mean a $\lambda_\Lambda$ which is an isomorphism up to prime-to-$p$ isogeny.
  \item $\eta^p$ is a level-$N^p$-structure, i.e. (if $S$ is connected), it is a $\pi_1(S,s)$-invariant $N^p$-orbit of symplectic similitudes $\eta^p\colon V_{\A_f^p}\to H_1(A_s,\A_f^p)$ (where $s$ is some geometric basepoint and $H_1(A_s,\A_f^p)$ with its $\pi_1(S,s)$-action corresponds to the Tate $\A_f^p$-module of $A$ (cf. \cite[\nopp 6.8]{rz}), which is a smooth $\A_f^p$-sheaf). Note that this forces the abelian schemes $A_\Lambda$ to be $(\dim_\Q V)$-dimensional.
  \end{enumerate}
\end{Remark}

\begin{Definition}
  Set $\Lambda^\S_{\Z_{(p)}}:=\Lambda^\S_{\Z_p}\cap V^\S_\Q=\prod_{i=-(r-1)-a}^{r-1}\Lambda_{\Z_{(p)}}^i$. We choose a lattice $\Lambda^\S_\Z\subseteq V^\S$ such that $\Lambda^\S_\Z\otimes_\Z\Z_{(p)}=\Lambda^\S_{\Z_{(p)}}$ and $\Lambda^\S_\Z\subseteq (\Lambda^\S_\Z)^\vee$.
\end{Definition}

\begin{Remark}\label{paragraph-moduli}
  Set $d:=\bigl|\left(\Lambda_\Z^\S\right)^\vee/\Lambda_\Z^\S\bigr|$.  By \cite[\nopp 2.3.3, 3.2.4]{kisin}, the integral model $\sS_J(\GSp(V^\S),S^{\S,\pm})$ is given by the moduli problem $(\Z_{(p)}\text{-schemes})\ni S\mapsto \left\{(A^\S,\lambda^\S,\epsilon^p)\right\}/\littlecong$, where
  \begin{enumerate}[(a)]
  \item $A^\S$ is an abelian scheme over $S$ up to $\Z_{(p)}$-isogeny,
  \item $\lambda^\S\colon A^\S\to \left(A^\S\right)^\vee$ is a polarization of degree $d$ (i.e., the polarization of the (well-defined) associated $p$-divisible group has degree $d$),
  \item $\epsilon^p$ is a level-$J^p$-structure, i.e. (if $S$ is connected), it is a $\pi_1(S,s)$-invariant $J^p$-orbit of symplectic similitudes $\epsilon^p\colon V^\S_{\A_f^p}\to H_1(A^\S_s,\A_f^p)$. Note that this forces the abelian schemes $A^\S$ to be $(\dim_\Q V^\S)$-dimensional.
  \end{enumerate}
\end{Remark}

This completes the descriptions of the moduli problems, and we turn to the question of the relationship between the two.  Consider (for appropriate $N^p,J^p$; see below) the morphism $\chi\colon\sS_N(\GSp(V),S^\pm) \to \sS_J(\GSp(V^\S),S^{\S,\pm})$ given on $S$-valued points by sending $(A,\bar\lambda,\eta^p)$ to $(A^\S,\lambda^\S,\epsilon^p)$, where
\begin{enumerate}[(a)]
\item $\displaystyle A^\S:=\prod_{i=-(r-1)-a}^{r-1}A_{\Lambda^i}$,
\item $\displaystyle \lambda^\S:=\prod_{i=-(r-1)-a}^{r-1}\left(\rho_{\left(\Lambda^i\right)^\vee,\Lambda^i}^\vee\circ \lambda_{\Lambda^i}\right)$,
\item $\epsilon^p$ is the product $\prod_{i=-(r-1)-a}^{r-1}\eta^p$, to be interpreted as the product over $\eta^p\colon V_{\A_f^p}\to H_1(A_{\Lambda^i,s},\A_f^p)\cong H_1(A_s,\A_f^p)$, where the isomorphism $H_1(A_{\Lambda^i,s},\A_f^p)\cong H_1(A_s,\A_f^p)$ is by definition the identity for some fixed $i=i_0$ and otherwise induced by the transition map $\rho_{\Lambda^{i},\Lambda^{i_0}}$. We need that $N^p$ is mapped into $J^p$ by $\GSp(V)\hookrightarrow\GSp(V^\S)$ for this to make sense.
\end{enumerate}

\begin{Lemma}\label{tate-faithful}
  Let $S$ be a scheme, $\ell\ne p$ prime numbers. If $\ell$ does not appear as a residue characteristic of $S$, then the Tate module functors
  \begin{align*}
    H_1(\_,\Z_\ell)&\colon (\text{abelian } S\text{-schemes})\to (\text{Ã©tale }\Z_\ell\text{-local systems on } S), \\
    H_1(\_,\Q_\ell)&\colon (\text{abelian } S\text{-schemes})\to (\text{Ã©tale }\Q_\ell\text{-local systems on } S)
  \end{align*}
  (cf. \cite[\nopp III, 5.4 and 6.2]{groth-bt} for precise definitions) are faithful.

  If only $p$ and $0$ appear as residue characteristics of $S$, then the Tate module functor
  \begin{equation*}
    H_1(\_,\A_f^p)\colon (\text{abelian } S\text{-schemes})\to (\text{Ã©tale }\A_f^p\text{-local systems on } S)
  \end{equation*}
  is faithful.
\end{Lemma}

\begin{Proof}
  First note that the statements about $H_1(\_,\Q_\ell)$ and $H_1(\_,\A_f^p)$ follows from the statement about $H_1(\_,\Z_\ell)$, which is why it is enough to only look at $H_1(\_,\Z_\ell)$.

  A homomorphism of abelian $S$-schemes $f\colon A\to B$ vanishes if and only if it vanishes over every (geometric) fiber of $S$: Indeed, if it vanishes fiberwise, then it is flat by the fiber criterion for flatness. Applying that criterion again we see that the closed immersion and fiberwise isomorphism $\ker(f)\hookrightarrow A$ is flat, which means that is an isomorphism.

  This way we are reduced to the case where $R$ is an (algebraically closed) field of characteristic different from $\ell$. In this setting the faithfulness is well-known (the salient point being that the $\ell$-primary torsion is dense).
\end{Proof}

\begin{Lemma}\label{dantzig}
  Let $H$ be a totally disconnected locally compact\footnote{By (our) definition, locally compact implies Hausdorff.} group (i.e., a locally profinite group) and let $N\subseteq H$ be a compact subgroup. Then
  \begin{equation*}
    N = \bigcap_{\substack{N\subseteq J \\J\subseteq H\text{ open compact subgrp.}}} J.
  \end{equation*}
\end{Lemma}

Note that this is (a variant of) a well-known theorem by van Dantzig if $N=\{1\}$ \cite{vandantzig}.

\begin{Proof}
  We make use of the following fact \cite[Prop.~3.1.7]{arhangel}: A Hausdorff space is locally compact and totally disconnected if and only if the open compact sets form a basis of the topology. (Van Dantzig's theorem is the group version of this, which talks only about a neighborhood basis of the identity and open compact \emph{subgroups}.)

  First we show that $N$ is contained in some open compact subset $K\subseteq H$. For every $x\in N$ choose a compact open neighborhood $x\in K_x\subseteq H$. This is possible by the fact cited above. Then there is a finite subset $I\subseteq N$ such that $N\subseteq \bigcup_{x\in I}K_x=:K$.

  Next, for every $x\in N$ choose an open neighborhood of the identity $U_x$ such that $xU_xK\subseteq K$.  With $N\subseteq U:=\bigcup_{x\in N}xU_x$ we obtain $UK\subseteq K$. Replacing $U$ by $U\cap U^{-1}$, we may moreover assume it is symmetric. The subgroup generated by $U$ is open (hence closed) and contained in $K$, hence is an open compact subgroup.

  Thus $N$ even is contained in an open compact sub\emph{group}; in other words, we may assume that $H$ is compact, i.e., is a profinite group.

  Then $H/N$ is compact\footnote{Hausdorff quotient spaces of compact spaces are compact again, but for ``locally compact'' the analogous statement is not true in general!} and totally disconnected\footnote{Take $x,y\in H$ such that $xN\ne yN$. We show that any subspace $S\subseteq H/N$ containing both $xN$ and $yN$ is disconnected. Let $U\subseteq H/N$ be a neighborhood of $xN$ not containing $yN$. Let $x\in V\subseteq \pi^{-1}(U)$ be open and compact, where $\pi\colon H\to H/N$ is the projection. Then $yN\notin \pi(V)\subseteq H/N$ is open and compact (hence closed) and we have $S=(\pi(V)\cap S)\sqcup S\setminus \pi(V)$ where both $\pi(V)\cap S$ and $S\setminus\pi(V)$ are open in $S$. This shows that $S$ is disconnected.} (i.e., is a Stone space). By the fact cited above,
  \begin{equation*}
    H/N \supseteq\{1\} = \bigcap_{L\subseteq H/N\text{ open compact subset}} L.
  \end{equation*}

  Observe that the quotient map $H\to H/N$ is proper to deduce
  \begin{equation*}
    N = \bigcap_{\substack{N\subseteq M\\M\subseteq H\text{ open compact subset}}} M.
  \end{equation*}

  Say $M$ is an open and compact subset of $H$ containing $N$.  As we have shown above, there is an open compact subgroup $J\subseteq H$ in between $N$ and $M$, and this is all we need to complete the proof.
\end{Proof}

\begin{Proposition}\label{diag-emb}
  For every compact open subgroup $N^p\subseteq \GSp(V)(\A_f^p)$
  \begin{equation*}
    \chi\colon\sS_N(\GSp(V),S^\pm) \to \sS_J(\GSp(V^\S),S^{\S,\pm})
  \end{equation*}
  is a well-defined morphism for all compact open subgroups $N^p\subseteq J^p\subseteq \GSp(V^\S)(\A_f^p)$ and is a closed immersion for all sufficiently small compact open subgroups $N^p\subseteq J^p\subseteq \GSp(V^\S)(\A_f^p)$.
\end{Proposition}

\begin{Proof}
  The fact that it's well-defined is clear from the construction.

  To show the second statement, as in \cite[Prop.~1.15]{travaux}, it is enough to show that
  \begin{equation*}
    \sS_{N_pN^p}(\GSp(V),S^\pm) \to \varprojlim_{J^p}\sS_{J_pJ^p}(\GSp(V^\S),S^{\S,\pm})
  \end{equation*}
  is a closed immersion, i.e., a proper monomorphism.

  We begin by proving that it is a monomorphism, i.e., injective on $S$-valued points ($S$ arbitrary $\Z_{(p)}$-scheme). So, say $(A_1,\lambda_1,\eta_1^p)$ and $(A_2,\lambda_2,\eta_2^p)$ both map to $(A^\S,\lambda^\S,\epsilon_{J^p}^p)$. That means precisely that there is an isomorphism of abelian $S$-schemes up to $\Z_{(p)}$-isogeny
  \begin{equation*}
    \phi\colon \prod_{i=-(r-1)-a}^{r-1}A_{1,\Lambda^i} \xrightarrow{\cong} \prod_{i=-(r-1)-a}^{r-1}A_{2,\Lambda^i}
  \end{equation*}
  such that
  \begin{equation*}
    \phi^\vee \circ \prod_{i=-(r-1)-a}^{r-1}\left(\rho_{2,\left(\Lambda^i\right)^\vee,\Lambda^i}^\vee\circ \lambda_{2,\Lambda^i}\right) \circ \phi = \prod_{i=-(r-1)-a}^{r-1}\left(\rho_{1,\left(\Lambda^i\right)^\vee,\Lambda^i}^\vee\circ \lambda_{1,\Lambda^i}\right)
  \end{equation*}
  and
  \begin{equation*}
    H_1(\phi,\A_f^p)\circ \epsilon_{1,J^p}^p = \epsilon_{2,J^p}^p   \mod{J^p}.
  \end{equation*}

  We claim that $\phi$ comes from isomorphisms
  \begin{equation*}
    \phi_i\colon A_{1,\Lambda^i} \xrightarrow{\cong} A_{2,\Lambda^i}.
  \end{equation*}

  Certainly there is but one candidate for $\phi_i$: define $\phi_i$ to be the composition
  \begin{equation*}
    A_{1,\Lambda^i}\xrightarrow{\mathrm{incl}} \prod_{i=-(r-1)-a}^{r-1}A_{1,\Lambda^i} \xrightarrow{\phi} \prod_{i=-(r-1)-a}^{r-1}A_{2,\Lambda^i} \xrightarrow{\mathrm{proj}} A_{2,\Lambda^i}.
  \end{equation*}

  Our claim then is that
  \begin{equation*}
    \phi = \prod_{i=-(r-1)-a}^{r-1}\phi_i.
  \end{equation*}

  Apply $H^1(\_,\A_f^p)$ on both sides. For the left hand side, we have
  \begin{equation*}
    H_1(\phi,\A_f^p) = \epsilon_{2,J^p}^p\circ \left(\epsilon_{1,J^p}^p\right)^{-1}   \mod{J^p}.
  \end{equation*}
  and the right hand side of this equation is block diagonal. So
  \begin{equation*}
    H_1(\phi,\A_f^p) = \prod_{i=-(r-1)-a}^{r-1}H_1(\phi_i,\A_f^p) \mod{J^p}.
  \end{equation*}

  Since (by Lemma~\ref{dantzig})
  \begin{equation*}
    N^p = \bigcap_{\substack{N_\ell\subseteq J_\ell\\ J_\ell\subseteq \GSp(V^\S)(\Q_\ell)\text{ cpt. open subgrp.}}}J_\ell,
  \end{equation*}
  it follows that (with $\ell\ne p$)
  \begin{equation*}
    H_1(\phi,\Q_\ell) = \prod_{i=-(r-1)-a}^{r-1}H_1(\phi_i,\Q_\ell) \mod{N_\ell},
  \end{equation*}
  hence (since $N_\ell$ acts block-diagonally) that $H_1(\phi,\Q_\ell)= \prod_{i=-(r-1)-a}^{r-1}H_1(\phi_i,\Q_\ell)$.

  Since $H_1(\_,\Q_\ell)$ is faithful (Lemma~\ref{tate-faithful}), this implies $\phi=\prod_{i=-(r-1)-a}^{r-1}\phi_i$, as desired.

  Next, consider the extension by zero of $\left(H_1(\rho_{1/2,\Lambda^j,\Lambda^i},\A_f^p)\right)_{i,j}$ (where for ``$1/2$'' either ``$1$'' or ``$2$'' can be plugged in) to a map $H_1(A^\S,\A_f^p)\to H_1(A^\S,\A_f^p)$. Under the isomorphism given by the $J^p$-level structure this corresponds, up to the $J^p$-action, to the map $V^\S_{\A_f^p}\to V^\S_{\A_f^p}$ given by mapping the $i$'th copy of $V_{\A_f^p}$ identically to the $j$'th copy and the rest to zero. Thus $\rho_{1/2,i,j}$ yield the same up to $J^p$ after applying $H_1(\_,\A_f^p)$, hence they are equal in the $\Z_{(p)}$-isogeny category.

  Consequently, $\chi$ is a monomorphism.

  For properness, we will use the valuative criterion. Let $R$ be a discrete valuation ring with field of fractions $K$ and assume that a $K$-point $A^\S=\prod_{i=-(r-1)-a}^{r-1} A_{\Lambda^i}$ with its additional structures coming from $(A_{\Lambda^i})_i$ extends to an $R$-point $\mathcal{A}^\S$. Consider the map $A^\S\to A_{\Lambda^{i_0}}\to A^\S$ where the first map is a projection and the second an inclusion. By the NÃ©ron mapping property, this extends to a map $\mathcal{A}^\S\to\mathcal{A}^\S$. Define $\mathcal{A}_{\Lambda^{i_0}}$ to be the image of this map.

  The NÃ©ron mapping property also allows us to extend the transition isogenies $\rho_{\Lambda^{i_0},\Lambda^{j_0}}\colon\allowbreak {A_{\Lambda^{j_0}}\to A_{\Lambda^{i_0}}}$, $i_0\leq j_0$,  the periodicity isomorphisms, and the polarization.

  Since $\pi_1(\Spec K)$ surjects onto $\pi_1(\Spec R)$ (see \stacks{0BQM}), extending the level structure away from $p$ is trivial.
\end{Proof}

\subsection{Local structure of the integral model}
\label{sec:local-structure}

\subsubsection{Generizations and irreducible components}
\label{sec:gener-irred-comp}

Let $\mathscr{X}\to \Spec \mathcal{O}_{\bE}$ be a flat
scheme locally of finite type; denote the special fiber by $X\to \Spec\bar\F_p$ and the generic fiber by $\mathcal{X}\to\Spec\bE$. We assume that $\mathcal{X}$ is locally integral (e.g. smooth).

For example, we can consider $(\mathscr{X},X,\mathcal{X})=(\mathscr{S}^-_K(G,X)_{\bOE}, \mathscr{S}^-_K(G,X)_{\bOE}\otimes_{\bOE}\bar\F_p,\allowbreak {\Sh_K(G,X)\otimes_E\bE})$.

Let $\bar x\in X(\bar\F_p)$. 

\begin{Lemma}\label{closed-generizations}
  There is a generization $x$ of $\bar x$ which lies in the generic fiber $\mathcal{X}$, and is a closed point in there, i.e., $x\in \mathcal{X}(L)$ for a finite extension $L/\bE$.
\end{Lemma}

\begin{Definition}\label{def-cpg}
  We shall call such a point $x$ a \defn{closed point generization} of $\bar x$ for short.
\end{Definition}

\begin{Proof}
  Due to flatness (going-down) there is \emph{some} generization in the generic fiber; call it $x_0$.

  By \stacks{053U} the following set is dense (and in particular non-empty) in the closure of $\{x_0\}$ in $\mathcal{X}$:
  \begin{equation*}
    \left\{ x\in\mathscr{X} \suchthat x\text{ is a specialization of } x_0\text{ and a closed point generization of } \bar x \right\}.
  \end{equation*}
\end{Proof}

\begin{Lemma}\label{closed-generizations-realization}
  Notation as in the preceding lemma.

  The specialization $x\leadsto \bar x$ can be realized by an ${\cal O}_L$-valued point of $\mathscr{X}$.
\end{Lemma}

\begin{Proof}
  First off, by \cite[\nopp 7.1.9]{EGA2}, it can be realized by a morphism $\Spec R=\{\eta,s\}\to\mathscr{X}$ of $\bOE$-schemes, where $R$ is a discrete valuation ring such that $L\cong\kappa(\eta)=\Quot(R)$ as field extensions of $\kappa(x)$.

  We hence get local homomorphisms of local rings $\bOE\to{\cal O}_{\mathscr{X},\bar x}\to R$.

  Thus the discrete valuation on $L$ defined by $R$ extends the discrete valuation on $\bE$. But there is but one such extension and its valuation ring is ${\cal O}_L$ (by definition).
\end{Proof}

\begin{Lemma}
  Mapping $x$ to the unique irreducible component of $\mathscr{X}$ that contains $x$ establishes a surjection from the set of closed point generizations $x$ of $\bar x$ to the set of irreducible components of $\mathscr{X}$ containing $\bar x$.
\end{Lemma}

\begin{Proof}
  If $x_0\in\mathcal{X}$ is a generization of $\bar x$, then $x_0$ lies in a unique irreducible component of $\mathscr{X}$ because $\mathcal{X}$ is locally irreducible. Hence the map described above is well-defined.

  Now for surjectivity: Given an irreducible component $C$ of $\mathscr{X}$ containing $\bar x$, let $x_0\in C$ be the generic point. Then $x_0$ must be in the generic fiber (else we would be able to find a generization in the generic fiber by going-down). Now go through the proof of Lemma~\ref{closed-generizations} with this particular choice of $x_0$.
\end{Proof}

\subsection{The local model}
\label{sec:local-model}

To give a very rough idea of what the \emph{local model} to be discussed in this section is supposed to accomplish: It should be an $\mathcal{O}_E$-scheme that is Ã©tale-locally isomorphic to $\sS_K(G,X)$, but easier to understand by virtue of being of a more ``linear-algebraic flavor''.  In actuality however, the theory of local models quickly gets quite complicated once one departs from the simplest examples.

\subsubsection{The Siegel case}
\label{sec:siegel-case-loc-mod}

We do start with the simplest example.

We consider the standard Iwahori subgroup $I_p\subseteq\GSp_{2g}(\Z_p)$, defined as the preimage of the standard Borel subgroup of $\GSp_{2g}(\F_p)$.  In terms of the building (cf. Remark~\ref{gsp-building}), it corresponds to the lattice chain $\mathcal{L}_\mathrm{full}$ given by
\begin{equation}\label{eq:siegel-full-lattice-chain}
  \begin{aligned}
    \Lambda^0=\Z_p^{2g}
    &\supsetneqq \Lambda^1=\Z_p^{2g-1}\oplus p\Z_p
    \supsetneqq \Lambda^2=\Z_p^{2g-2}\oplus p\Z_p^2 \\
    &\supsetneqq \dotsb
    \supsetneqq \Lambda^{2g-1} = \Z_p\oplus p\Z_p^{2g-1}
    \supsetneqq p\Lambda^0=p\Z_p^{2g}
  \end{aligned}
\end{equation}
of period length $2g$.

Consider a subset $J=\{j_0>\dotsb>j_{m-1}\}\subseteq \{1,\dotsc,2g\}$ such that for each $j\in J$ with $1\leq j\leq 2g-1$ also $2g-j\in J$, and let $K_p$ be the parahoric subgroup associated with the partial lattice chain $\mathcal{L}\subseteq\mathcal{L}_\mathrm{full}$ obtained from $\left\{\Lambda^j \suchthat j\in J\right\}$. 

Define a scheme $\tsS_K(G,X)$ over $\sS_K(G,X)$ as follows:
\begin{equation*}
  \tsS_K(G,X)(S) =
\left\lbrace (A,\bar\lambda,\eta^p,\tau) \;\middle|\;
  \begin{tabular}{@{}l@{}}
    $(A,\bar\lambda,\eta^p)\in\sS_K(\GSp_{2g},S^\pm)(S)$, \\
    $\tau\colon H_\mathrm{dR}^1(A)\xrightarrow{\sim} \mathcal{L}\otimes\mathcal{O}_S$ isomorphism of lattice chains
   \end{tabular}
  \right\rbrace
\end{equation*}
for every $\Z_p$-scheme $S$.

By \cite[Appendix to Chap.~3]{rz}, $\tsS_K(G,X)\to \sS_K(G,X)$ is a Zariski torsor under the automorphism group of $\mathcal{L}$, i.e., the Iwahori group scheme.

This motivates the definition of the local model $M^\mathrm{loc}_{K_p}\to\Spec\Z_p$ as the ``moduli space of Hodge filtrations''; more precisely:

\begin{Remark}\label{loc-mod-siegel} \textnormal{(See \cite[91]{goertz-sympl}.)}
  $M^\mathrm{loc}_{K_p}(S)$ is the set of isomorphism classes of commutative diagrams
  \begin{equation*}
    \begin{tikzpicture}[node distance=4cm, auto]
      \node (X) at (0,2) {$\Lambda^{j_0}_S$};
      \node (Xs) at (0,0) {$\mathcal{F}^{j_0}$};
      \node (Z1) at (2.5,2) {$\Lambda^{j_1}_S$};
      \node (Zs1) at (2.5,0) {$\mathcal{F}^{j_1}$};
      \node (Z2) at (5,2) {$\dotsb$};
      \node (Zs2) at (5,0) {$\dotsb$};
      \node (Z) at (7.5,2) {$\Lambda^{j_0}_S$};
      \node (Zs) at (7.5,0) {$\mathcal{F}^{j_0}$};
      \node (O) at (10,2) {$\Lambda^{j_{m-1}}_S$};
      \node (Os) at (10,0) {$\mathcal{F}^{j_{m-1}}$};
      
      \draw[->] (X) to (Z1);
      \draw[->] (Z1) to (Z2);
      \draw[->] (Z2) to (Z);
      \draw[->] (Z) to node{$\cdot p$} (O);
      \draw[<-left hook] (Z) to (Zs);
      \draw[<-] (Zs) to  (Zs2);
      \draw[<-] (Zs2) to (Zs1);
      \draw[<-] (Zs1) to (Xs);
      \draw[<-] (Os) to (Zs);
      \draw[left hook->] (Xs) to (X);
      \draw[left hook->] (Zs1) to (Z1);
      \draw[left hook->] (Os) to (O);
    \end{tikzpicture}
  \end{equation*}
  with $\Lambda^{j}_S:=\Lambda^j\otimes_{\Z_p}\mathcal{O}_S$, $\mathcal{F}^j\subseteq\Lambda^j_S$ locally direct summand
  of rank $g$, such that for all $j\in J$, $\mathcal{F}^j\to \Lambda^j_S\overset{\psi}\cong (\Lambda^{2g-j}_S)^* \to (\mathcal{F}^{2g-j})^*$ vanishes, $\psi$ being the symplectic pairing.
\end{Remark}

By Grothendieck-Messing theory, one obtains a diagram
\begin{equation*}
  \begin{tikzpicture}[node distance=4cm, auto]
    \node (X) at (2.5,3) {$\tsS_K(G,X)$};
    \node (Z) at (0,0) {$\sS_K(G,X)$};
    \node (Xs) at (5,0) {$M^\mathrm{loc}_K$};
    
    \draw[->] (X) to node {smooth of rel. dim. $\dim \Aut(\mathcal{L})$} (Xs);
    \draw[->] (X) to node[swap] {$\Aut(\mathcal{L})$-torsor} (Z);
  \end{tikzpicture}
\end{equation*}

Since both morphisms in this diagram are smooth of the same dimension, it follows that for every finite field extension $\F_q/\F_p$ and every point $x\in\sS_K(G,X)(\F_q)$, there exists a point $y\in M^\mathrm{loc}_K(\F_q)$ and an isomorphism $\mathcal{O}_{\sS_K(G,X),x}^h\cong\mathcal{O}_{M^\mathrm{loc}_K,y}^h$ of henselizations.

In many (P)EL situations one has similar descriptions with the obvious extra structures. Sometimes however the so-called ``naive'' local models so obtained additionally need to be flattened, which leaves one without any self-evident moduli interpretation.

\subsubsection{The relation between the integral and the local model}
\label{sec:int-model-loc-model}

Generalizing the Siegel example, we axiomatically characterize the relationship between the integral model of the Shimura variety and its local model: One wants a \defn{local model diagram}, i.e., a diagram of $\mathcal{O}_E$-schemes functorial in $K$
\begin{equation}\label{eq:loc-mod-diagr}
  \begin{tikzpicture}[node distance=4cm, auto]
    \node (X) at (2.5,3) {$\tsS_K(G,X)$};
    \node (Z) at (0,0) {$\sS_K(G,X)$};
    \node (Xs) at (5,0) {$M^\mathrm{loc}_K$};
    
    \draw[->] (X) to node {equivariant and smooth of rel. dim. $\dim \mathcal{G}_{\mathcal{O}_E}$} (Xs);
    \draw[->] (X) to node[swap] {$\mathcal{G}_{\mathcal{O}_E}$-torsor} (Z);
  \end{tikzpicture}
\end{equation}
where $M^\mathrm{loc}_K$ is a projective flat $\mathcal{O}_E$-scheme with an action of $\mathcal{G}\otimes_{\Z_p}\mathcal{O}_E$ and generic fiber the canonical model of $G_{\bar\Q_p}/P_{\mu^{-1}}$ over $E$.

By Kisin-Pappas \cite{kisin-pappas} we do actually have such a diagram in our situation.

\subsubsection{The Pappas-Zhu construction}
\label{sec:papp-zhu-constr}

In \cite{pappas-zhu}, Pappas and Zhu give a construction of the local model in quite a general context, in particular with no assumptions going beyond our running assumptions \ref{std-assum}.

\begin{Remark}
  To this end, they construct an affine smooth group scheme $\ucG_K\to \A^1_{\Z_p}=\Spec \Z_p[t]$ with the following key properties:
  \begin{enumerate}[(1)]
  \item $\ucG_K$ has connected fibers,
  \item $\ucG_K$ is reductive over $\Spec \Z_p[t^{\pm 1}]$,
  \item $\ucG_K\otimes_{\Z_p[t],t\mapsto p}\Z_p\cong \cG_K$, in particular
    \begin{itemize}
    \item $\ucG_K\otimes_{\Z_p[t],t\mapsto p}\Q_p\cong G_{\Q_p}$ and
    \item $\ucG_K\otimes_{\Z_p[t]}\F_p:=\ucG_K\otimes_{\Z_p[t],t\mapsto 0}\F_p\cong \cG_K\otimes\F_p$,
    \end{itemize}
  \item  $\ucG_K\otimes_{\Z_p[t]}\Q_p\llbracket t\rrbracket$ is parahoric for $\ucG_K\otimes_{\Z_p[t]}\Q_p\llparen t\rrparen$,
  \item $\ucG_K\otimes_{\Z_p[t]}\F_p\llbracket t\rrbracket$ is parahoric for $\ucG_K\otimes_{\Z_p[t]}\F_p\llparen t\rrparen$.
  \end{enumerate}
\end{Remark}

\begin{DefinitionRemark}\label{pappzhu}
  Let $X_\mu$ be the canonical model of $G_{\bar\Q_p}/P_{\mu^{-1}}$ over $E$, where for a cocharacter $\nu$ one defines $P_{\nu}:=\{ g \in G \suchthat \lim_{t\to 0}\nu(t)g\nu(t)^{-1}\text{ exists}\}$.

  Let $S_\mu$ be the closed subvariety of $\Gr_{G}\times_{\Q_p}E$ with
  \begin{equation*}
    S_\mu(\bar\Q_p) = G(\bar\Q_p\llbracket t\rrbracket)\mu(t)G(\bar\Q_p\llbracket t\rrbracket) / G(\bar\Q_p\llbracket t\rrbracket).
\end{equation*}

  Then $S_\mu$ can be $G_E$-equivariantly identified with $X_\mu$.
\end{DefinitionRemark}

\begin{Definition}
  The local model $M^\mathrm{loc}_{G,\mu,K}$ now is defined to be the Zariski closure of $X_\mu\subseteq \Gr_{G}\times_{\Q_p} E$ in $\Gr_{\ucG_K,\Z_p}\otimes_{\Z_p}\mathcal{O}_E$, where $\Gr_{\ucG_K,\Z_p}:=\Gr_{\ucG_K,\A_{\Z_p}^1}\otimes_{\A^1_{\Z_p},\,u\mapsto p}\Z_p$ is a base change of the global affine Gra\ss{}mannian as defined in \cite{pappas-zhu}.
\end{Definition}

\section{EKOR strata and zips in the case of parahoric reduction}
\label{part:ekor-strata-zips}

\begin{Notation}
  We still fix a Shimura datum $(G,X)$ of Hodge type, a parahoric subgroup $K_p\subseteq G(\Q_p)$ (associated with a Bruhat-Tits group scheme $\mathcal{G}=\cG_K=\cG_{K_p}\to\Spec\Z_p$ associated with a facet $\ff$) and a sufficiently small open compact subgroup $K^p\subseteq G(\A_f^p)$. Define $\ocG_K:=\cG_K\otimes_{\Z_p}\kappa$.

  We also keep up our standard assumptions~\ref{std-assum}.
\end{Notation}

We now want to discuss the EKOR stratification on the special fiber of the integral model with parahoric level structure. The EKOR stratification interpolates between the Ekedahl-Oort (EO) and the Kottwitz-Rapoport (KR) stratification (see Remark~\ref{ekor-interpoliert} below for a precise formulation). We begin by explaining the basics about these stratifications and the combinatorics involved in the first section of this chapter.

\subsection{The Ekedahl-Oort, Kottwitz-Rapoport and EKOR stratifications}
\label{sec:eo-kr-ekor}

\subsubsection{Iwahori-Weyl group and the admissible subset}
\label{sec:iwahori-weyl-group}

\begin{Notation}\label{iwahori-notation}
  \begin{enumerate}[(1)]
  \item We fix an Iwahori subgroup $I_p\subseteq K_p$, i.e., $I_p$ is the group of $\Z_p$-points of the parahoric group scheme $\mathcal{I}$ associated with an alcove $\fa$ (facet of maximal dimension) such that $\ff\subseteq\overline{\fa}$. As usual, we also define $\bI:=\mathcal{I}(\bZ_p)\subseteq \bK$.
  \item Let $T\subseteq G$ be a maximal torus such that $T_{\bQ_p}$ is contained in a Borel subgroup of $G_{\bQ_p}$\footnote{Note that by Steinberg's theorem, $G_{\bQ_p}$ is quasi-split. \cite[Chap.~III, Â§\,2]{serre-galois}} and let $S$ be the maximal $\bQ_p$-split torus contained in $T_{\bQ_p}$. We can and do choose $T$ such that the alcove $\fa$ is contained in the apartment associated with $S$. By $N$ we denote the normalizer of $T$.
  \item Let $(V,R)$ be the relative root system of $(G_{\bQ_p},T_{\bQ_p})$, i.e., $V$ is the $\R$-vector space $X^*_{\bQ_p}(T_{\bQ_p})\otimes_\Z\R$ and $R\subseteq X^*_{\bQ_p}(T_{\bQ_p})$ is (as usual) such that we have a decomposition
    \begin{equation*}
      \fg:=\Lie(G_{\bar\Q_p})=\Lie(T_{\bar\Q_p})\oplus\bigoplus_{\alpha\in R}\fg_\alpha.
    \end{equation*}
    Contrary to the absolute situation, $\dim\fg_\alpha$ may be greater than $1$.
  \end{enumerate}
\end{Notation}

\begin{Definition}\label{weyl}
  \begin{enumerate}[(1)]
  \item \label{item:weyl-relweyl} The \defn{(finite relative) Weyl group} of $G$ (over $\bQ_p$) is $W:=N(\bQ_p)/T(\bQ_p)$. It is the Weyl group of the root system $(V,R)$, i.e., the group generated by the orthogonal reflections through the hyperplanes defined by the elements of $R$.
  \item As described in \cite[\nopp 1.2.3]{landvogt}, one defines a set of affine roots $R_\mathrm{aff}\supseteq R$ on $V$ using the valuation on $\bQ_p$. By $W_a\subseteq\mathrm{Aff}(V^*)=\GL(V^*)\ltimes V^*$ we denote the \defn{affine Weyl group} of the affine root system $(V,R_\mathrm{aff})$, i.e., the group generated by the orthogonal reflections through the affine hyperplanes defined by the elements of $R_\mathrm{aff}$.
  \item $\tW := N(\bQ_p)/(T(\bQ_p)\cap \bI)$ is the \defn{Iwahori-Weyl group}.
  \item $W_K:= (N(\bQ_p)\cap \bK)/(T(\bQ_p)\cap \bI)\subseteq\tW$. (Recall that $\bK=\mathcal{G}(\bZ_p)$.)
  \end{enumerate}
\end{Definition}

\begin{Remarks}
  \begin{enumerate}[(1)]
  \item We have $W\subseteq W_a$. With the systems of generators indicated above, $W$ and $W_a$ become (affine) Coxeter groups; in particular we can talk about reduced words and have length functions, cf. \cite{brenti}.
  \item $W_I$ is the trivial group.
  \end{enumerate}
\end{Remarks}

\begin{Proposition}\textnormal{\cite[Prop.~8]{hainesrap}}\label{bt-decomp}
  The Bruhat-Tits decomposition
  \begin{equation*}
    G(\bQ_p) = \bigcup_{w\in \tW}\bK w\bK
  \end{equation*}
  identifies
  \begin{equation*}
    \bK\backslash G(\bQ_p) / \bK \cong W_K\backslash \tW/W_K.
  \end{equation*}
\end{Proposition}

\begin{Proposition}\textnormal{\cite[Prop.~13]{hainesrap}}
  Let $\bK$ be the maximal parahoric subgroup of $G(\bQ_p)$ associated with a special vertex in the apartment corresponding to $S$. Then $W_K\to W$ is an isomorphism and $\tW\cong W\ltimes X_*(T)_{\Gal(\bar\Q_p/\bQ_p)}$.\footnote{Notation: Let $\Gamma$ be a group and $M$ a $\Z[\Gamma]$-module. Then $M_\Gamma:=\Z\otimes_{\Z[\Gamma]}M=M/\langle \gamma m-m \suchthat \gamma\in\Gamma,\; m\in M\rangle$ is the module of $\Gamma$-coinvariants of $M$.}
\end{Proposition}

\begin{Notation}
  We denote the map $X_*(T)_{\Gal(\bar\Q_p/\bQ_p)}\to \tW$ of the proposition by $\nu\mapsto t_\nu$.
\end{Notation}

\begin{Proposition}\textnormal{\cite[Lemma~14]{hainesrap}}
  Let $\Omega\subseteq\tW$ be the subgroup consisting of those elements that preserve the base alcove $\fa$.
  
  There is an exact sequence
  \begin{equation*}
    1 \to W_a \to \tW \to \Omega \to 1,
  \end{equation*}
  with a canonical right splitting (namely the inclusion $\Omega\hookrightarrow \tW$), i.e., $\tW \cong W_a\rtimes\Omega$.
\end{Proposition}

\begin{Definition}
  The semidirect product decomposition of the preceding proposition means that $\tW$ is a ``quasi-Coxeter'' group. In practical terms, this means:
  
  \begin{enumerate}[(1)]
  \item We define a length function $\ell$ on $\tW$ as follows: $\ell(w_a,\omega):=\ell(w_a)$ for all $w_a\in W_a$ and $\omega\in\Omega$, where on the right hand side we use the length function of the affine Coxeter group $W_a$.

    Note that $\Omega=\ell^{-1}(0)$.
  \item Likewise, we extend the Bruhat partial order from $W_a$ to $\tW$ by defining
    \begin{equation*}
      (w_{a,1},\omega_1) \leq (w_{a,2},\omega_2) \coloniff w_{a,1}\leq w_{a,2}\text{ and } \omega_1=\omega_2.
    \end{equation*}

    Note that $w_1\leq w_2$ ($w_1,w_2\in\tW$) implies $\ell(w_1)\leq\ell(w_2)$.
  \end{enumerate}
\end{Definition}

\begin{Definition}
  \begin{enumerate}[(1)]
  \item Let $\{\mu\}$ be a $W_\mathrm{abs}$-conjugacy class  of geometric cocharacters of $T$ (cf. Remark~\ref{conj-class-cochar}), $W_\mathrm{abs}:=N(\bar\Q_p)/T(\bar\Q_p)$ being the absolute Weyl group. Let $\bar\mu\in X_*(T)_{\Gal(\bar\Q_p/\bQ_p)}$ be the image of a cocharacter in $\{\mu\}$ whose image in $X_*(T)\otimes_\Z\R$ is contained in the closed Weyl chamber corresponding to some Borel subgroup of $G$ containing $T$ and defined over $\bQ_p$.
  \item $\Adm(\mu):=\Adm(\{\mu\}):=\{w\in\tW\suchthat w \leq q t_{\bar\mu}q^{-1}=t_{q\bar\mu}\text{ for some } q\in W\}$ is the \defn{$\{\mu\}$-admissible subset} of $\tW$.
  \item $\Adm(\{\mu\})^K:= W_K\Adm(\{\mu\})W_K\subseteq \tW$.
  \item $\Adm(\{\mu\})_K:=\KR(K,\{\mu\}):=W_K\backslash \Adm(\{\mu\})^K/W_K\subseteq W_K\backslash \tW/W_K$.
  \item Define $^K\tW\subseteq\tW$ to be the set of representatives of minimal length for the quotient $W_K\backslash \tW$.
  \item $^K\Adm(\{\mu\}):=\EKOR(K,\{\mu\}):=\Adm(\{\mu\})^K\cap {}^K\tW\subseteq {}^K\tW$.
  \end{enumerate}
\end{Definition}

\begin{Lemma}\textnormal{(See \cite[Thm.~1.2.2]{SYZnew}.)} 
  $^K\Adm(\{\mu\})=\Adm(\{\mu\})\cap{}^K\tW$.
\end{Lemma}

\subsubsection{Kottwitz-Rapoport stratification}
\label{sec:kottw-rapop-strat}

Recall from Section~\ref{sec:int-model-loc-model} that we have an integral model and a local model diagram
\begin{equation*}
  \sS_K \leftarrow \tsS_K \to M^\mathrm{loc}_K
\end{equation*}
or, equivalently, a (smooth) morphism of stacks $\sS_K \to [\cG_K \backslash M^\mathrm{loc}_K]$ (over $\mathcal{O}_E$).

As explained in Section~\ref{sec:papp-zhu-constr}, by the construction in \cite{pappas-zhu}, the special fiber $M^\mathrm{loc}_K\otimes\kappa$ of $M^\mathrm{loc}_K$ is a closed subscheme of the affine flag variety $\Gr_{\ucG_K\otimes_{\Z_p}\kappa}=\mathcal{F}l_{\ucG_K\otimes\kappa\llbracket t\rrbracket}$, which is the ind-projective ind-scheme over $\kappa$ given as the fpqc sheafification (which exists in this case!) of the presheaf $R\mapsto \ucG_K(R\llparen t\rrparen)/\ucG_K(R\llbracket t\rrbracket)$.

\begin{Definition}
  Define $L^+(\ucG_K\otimes\kappa\llbracket t\rrbracket)$ to be the $\kappa$-functor sending a $\kappa$-algebra $R$ to $\ucG_K(R\llbracket t\rrbracket)$.

  We let $L^+(\ucG_K\otimes\kappa\llbracket t\rrbracket)$ act on $\Gr_{\ucG_K\otimes_{\Z_p}\kappa}$ from the left and call this action $a$ (within this subsection).  The orbits of this action on $\Gr_{\ucG_K\otimes_{\Z_p}\bar\kappa}$ are the \defn{Schubert cells}.
\end{Definition}

\begin{Remarks}\label{schubert-cell-indexing}
  \begin{enumerate}[(1)]
  \item \label{item:schubert-cell-indexing}The Schubert cells can be indexed by $W_K\backslash \tW/W_K$ by Proposition~\ref{bt-decomp} with the following in mind: Strictly speaking, using the Bruhat-Tits decomposition here, we arrive at something involving the Iwahori-Weyl group of $\ucG_K\otimes\bar{\kappa}\llparen t\rrparen$. However, by \cite[\nopp 9.2.2]{pappas-zhu}, this is isomorphic to the Iwahori-Weyl group of $G_{\bQ_p}$.
  \item $M^\mathrm{loc}_K\otimes\bar{\kappa}$ is a union of Schubert cells, namely of those indexed by $\KR(K,\{\mu\}):=W_K\backslash (W_K\Adm(\{\mu\})W_K)/W_K$, cf. \cite[Theorem~9.3]{pappas-zhu}.
  \end{enumerate}
\end{Remarks}

\begin{Remark}
  By construction, $M^\mathrm{loc}_K$ has an action $b$ of $\ucG_K\otimes_{\Z_p[t],t\mapsto p}\mathcal{O}_E\cong \cG_K\otimes\mathcal{O}_E$.

  For $w\in\tW$ choose a representative $\dot{w}\in L\ucG_K(\bar\kappa)$ (with Remark~\ref{schubert-cell-indexing}\,\eqref{item:schubert-cell-indexing} in mind) and let $e_0\in\Gr_{\ucG_K\otimes_{\Z_p}\kappa}$ be the distinguished base point (associated with the identity). For $w \in W_K\Adm(\{\mu\})W_K$, the orbit map of $\dot{w}\cdot e_0$ for the action $a$ factors through the homomorphism $L^+(\ucG_K\otimes\bar\kappa\llbracket t\rrbracket) \to \cG_K\otimes\kappa \cong \ucG_K\otimes\bar\kappa$.

  The orbits associated with the two $\cG_K\otimes\kappa$-actions $a$ and $b$ on
  $M^\mathrm{loc}_K\otimes\kappa$
  agree. The orbits of the $\cG_K\otimes\kappa$-action on $M^\mathrm{loc}_K\otimes\kappa$ are indexed by $\KR(K,\{\mu\})$.
\end{Remark}

\begin{Definition}
  The stratifications thus obtained on $M^\mathrm{loc}_K\otimes\kappa$ and $\sS_K\otimes\kappa$ are called \defn{Kottwitz-Rapoport stratifications}. That is to say that  Kottwitz-Rapoport strata on $\sS_K\otimes\kappa$ are by definition pullbacks of Kottwitz-Rapoport strata on $M^\mathrm{loc}_K$, which in turn are $\cG_K\otimes\kappa$-orbits.
\end{Definition}

\subsubsection{Ekedahl-Oort stratification}
\label{sec:eked-oort-strat}

The Ekedahl-Oort stratification is only defined in the case of good reduction, i.e., if $K_p$ is hyperspecial or, equivalently, if $\cG_K$ is a \emph{reductive} model of $G_{\Q_p}$. Then $G_{\Q_p}$ splits over $\bQ_p$ (by definition of ``hyperspecial'', cf. \cite[\nopp 1.10.2]{tits}).

We therefore put ourselves in the situation of good reduction for this subsection.

\begin{Remark}
  Then $W$ as defined in Definition~\ref{weyl}\,\eqref{item:weyl-relweyl} agrees with the absolute Weyl group of $G_{\Q_p}=\cG_K\otimes\Q_p$, which in turn agrees with the absolute Weyl group of $\ocG_K:=\cG_K\otimes\kappa$, cf. \cite[App.~A.5]{vieh-wed}.
\end{Remark}

\begin{Definition}
  Define $I$ to be the type (interpreted as a subset of simple reflections) of the parabolic subgroup of $G_{\Q_p}$ defined by $\mu^{-1}$ (cf. Remark~\ref{pappzhu}), and $^IW\subseteq W$ to be the system of representatives of the quotient group $W_I\backslash W$ containing the element of least length of every coset.
\end{Definition}

\begin{Theorem} \textnormal{\cite{mw,PWZ,PWZ-AZD,zhangEO}}
  There is a smooth algebraic stack $\ocG_K\Zip_\kappa:=\ocG_K\Zip^\mu_\kappa$ over $\kappa$ with underlying topological space $^IW$ together with a smooth morphism
  \begin{equation*}
    \sS_K\otimes\kappa \to \ocG_K\Zip^\mu_\kappa.
  \end{equation*}

  The stratification of $\sS_K\otimes\kappa$ thus obtained is the \defn{Ekedahl-Oort stratification}.
\end{Theorem}

\subsubsection{EKOR stratification}
\label{sec:ekor-stratification}

\begin{Definition}
  Let $L$ be a valued field extension of $\Q_p$ with ring of integers $\mathcal{O}$, maximal ideal $\fm$ and residue field $\lambda$.

  The \defn{pro-unipotent radical} of $\cG_K(\mathcal{O})$ is
  \begin{equation*}
    \cG_K(\mathcal{O})_1 := \{ g \in \cG_K(\mathcal{O}) \suchthat (g \mod \fm)\in \bar{R}_K(\lambda) \},
  \end{equation*}
  where $\bar{R}_K$ is the unipotent radical of $\cG_K\otimes_{\Z_p}\lambda$.
\end{Definition}

In particular, if $K$ is hyperspecial, then $\cG_K(\mathcal{O})_1=\ker(\cG_K(\mathcal{O})\to \cG_K(\lambda))$.

Also, $\obK:=\bK/\bK_1\cong \ocG_K^\mathrm{rdt}(\bar\F_p)$, where $\ocG_K^\mathrm{rdt}$ is the maximal reductive quotient of $\ocG_K:=\cG_K\otimes\kappa$.

\begin{Remark} \textnormal{\cite[after~Cor.~6.2]{he-rapo}}\label{iw-diagr}
  We have a commutative diagram
  \begin{equation*}
    \begin{tikzpicture}[node distance=4cm, auto]
      \node (X) at (0,2) {$G(\bQ_p)/\bK_\sigma(\bK_1\times\bK_1)$};
      \node (Z) at (5,2) {$^K\tW$};
      \node (Zs) at (5,0) {$W_K\backslash \tW/W_K.$};
      \node (Xs) at (0,0) {$\bK\backslash G(\bQ_p)/\bK$};
      
      \draw[->] (X) to (Xs);
      \draw[->] (Z) to  (Zs);
      \draw[->] (X) to  (Z);
      \draw[->] (Xs) to  (Zs);
    \end{tikzpicture}
  \end{equation*}
\end{Remark}

Consider the map
\begin{equation*}
  v_K\colon \sS_K\otimes\kappa\to G(\bQ_p)/\bK_\sigma(\bK_1\times\bK_1),
\end{equation*}
which is the composition of the central leaves map $\Upsilon_K\colon \sS_K\otimes\kappa\to G(\bQ_p)/\bK_\sigma$ (see \cite{article-leaves}) with the projection $G(\bQ_p)/\bK_\sigma\to G(\bQ_p)/\bK_\sigma(\bK_1\times\bK_1)$. The Kottwitz-Rapoport map $\lambda_K\colon {\sS_K\otimes\kappa} \to \bK\backslash G(\bQ_p)/\bK$ factors through this map.

The fibers of $v_K$ are called \defn{EKOR strata}. By \cite[Thm.~6.15]{he-rapo}, they are locally closed subsets of $\sS_K\otimes\kappa$.

\begin{Remarks}
  \begin{enumerate}[(1)]
  \item One can explicitly express the image of a EKOR stratum under a change-of-parahoric map as a union of EKOR strata on the target \cite[Prop.~6.11]{he-rapo}.
  \item The closure of an EKOR stratum is a union of EKOR strata and one can explicitly describe the associated order relation \cite[Thm.~6.15]{he-rapo}.
  \end{enumerate}
\end{Remarks}

\begin{Remark}\label{ekor-interpoliert}
  In the hyperspecial case, the EKOR stratification agrees with the Ekedahl-Oort stratification. In the Iwahori case, it agrees with the Kottwitz-Rapoport stratification ($^K\tW=\tW=W_K\backslash \tW/W_K$ in that case).
\end{Remark}

By definition, the EKOR stratification always is a refinement of the Kottwitz-Rapoport stratification.  So one way of approaching the EKOR stratification is by looking at a fixed Kottwitz-Rapoport stratum and trying to understand how it is subdivided into EKOR strata.

To get this started, let us recall some calculations from the proof of \cite[Thm.~6.1]{he-rapo}.

Fixing a Kottwitz-Rapoport stratum means restricting our view to $\bK w\bK/\bK_\sigma$ rather than the whole of $G(\bQ_p)/\bK_\sigma$, for some fixed $w\in\KR(K,\{\mu\})$. The EKOR strata in the Kottwitz-Rapoport stratum associated with $w$ are therefore indexed by $\bK w\bK/\bK_\sigma(\bK_1\times\bK_1)$.

Define $\sigma':=\sigma\circ\Ad(w)$ and consider the bijection
\begin{align*}
  \bK/(\bK\cap w^{-1}\bK w)_{\sigma'} &\xrightarrow{\sim} \bK w\bK/\bK_\sigma, \\
  k &\mapsto wk, \\
  k_2\sigma(k_1) &\mapsfrom k_1wk_2.
\end{align*}

Let $J$ be the set of simple affine reflections in $W_K$, let $\bar B$ be the image of $\bI$ in $\obK$ and $\bar T\subseteq \bar B$ the maximal torus. Set $J_1:=J\cap w^{-1}Jw$.

\begin{Proposition}\label{propmorris}
  \textnormal{(See \cite[Lemma~3.19]{morris}.)}
  The image of $\bK \cap w^{-1} \bK w$ in $\obK$ is $\bar{P}_{J_1}$, i.e., the standard parabolic subgroup of $\obK$ associated with $J_1$.
\end{Proposition}

\begin{Remark}
  He and Rapoport invoke Carter's book \cite{carter} at this point, which primarily pertains to the case of (usual) BN-pairs attached to reductive groups. Morris \cite{morris} shows that the relevant results carry over likewise to the case of generalized (or affine) BN-pairs.
\end{Remark}

Then we get a map
\begin{align*}
  \bK w\bK/\bK_\sigma \to \bK/(\bK\cap w^{-1}\bK w)_{\sigma'}  &\to \obK / (\bar{P}_{J_1})_{\sigma'} \\
  &\to \obK / (\bar{L}_{J_1})_{\sigma'}(\bar{U}_{J_1})_{\sigma'} \to \obK / (\bar{L}_{J_1})_{\sigma'}(\bar{U}_{J_1}\times \bar{U}_{\sigma'(J_1)}),
\end{align*}
which factors through a bijection
\begin{equation*}
  \bK w\bK/\bK_\sigma(\bK_1\times\bK_1) \xrightarrow{\sim}
  \obK / (\bar{L}_{J_1})_{\sigma'}(\bar{U}_{J_1}\times \bar{U}_{\sigma'(J_1)}) \cong \ocG_K^\mathrm{rdt}\Zip^{\mathcal{Z}_w}(\bar\F_p)/{\littlecong}.
\end{equation*}

Here, $\mathcal{Z}_w$ is the (connected) algebraic zip datum $\mathcal{Z}_w=(\ocG^\mathrm{rdt},\bar{P}_{J_1},\bar{P}_{\sigma'(J_1)},\sigma')$, as described in \cite{SYZnew}. In \cite{SYZnew}, Shen, Yu and Zhang show that this observation ``globalizes'' (with the drawback that ``global'' here still just refers to the Kottwitz-Rapoport stratum\footnote{They also give another ``globalization''; the drawback there being that it only works after perfection.}) in a pleasant way. To wit, one gets a smooth morphism \cite[Theorem~A]{SYZnew}
\begin{equation*}
  \zeta_w\colon \osS_K^w \to \ocG_K^\mathrm{rdt}\Zip^{\mathcal{Z}_w}_\kappa
\end{equation*}
(the source being a Kottwitz-Rapoport stratum).

\subsection{\texorpdfstring{$\ocG_K$}{GK}-zips in the Siegel case}
\label{sec:G_K-zips}

Here we work with the Siegel Shimura datum, cf. Example~\ref{sympl-ex}.

\subsubsection{Preliminaries}

\begin{Notation}
  Fix $p\ne 2$,\footnote{As in \cite{rz}, the principal reason for this restriction is our use of the equivalence between alternating and skew-symmetric. See Definition~\ref{siegel-chain}\,\eqref{item:siegel-chain-pol}.} $g\in\Z_{\geq 1}$ and a subset $J\subseteq \Z$ with $J=-J$ and $J+2g\Z=J$. Associated with $J$ is the partial lattice chain $\left\{\Lambda^j \suchthat j\in J\right\}$, where $\Lambda^j$ are defined as in equation~\eqref{eq:siegel-full-lattice-chain}. Let $K_p$ be the corresponding parahoric subgroup of $\GSp_{2g}(\Q_p)$, i.e., the stabilizer of said lattice chain. It contains the Iwahori subgroup $I_p$ associated with the full lattice chain~\eqref{eq:siegel-full-lattice-chain}. For the maximal torus $T$ we take the usual diagonal (split) torus.
\end{Notation}

\begin{Remark}\label{root-stuff-gsp}
  The Weyl group is
  \begin{align*}
    W&=\{\pi\in S_{2g}=\Aut(\{\pm 1,\pm 2,\dotsc,\pm g\}) \suchthat \pi(-n)=-\pi(n)\text{ for } n=1,2,\dotsc,g\}\\
     &\cong S_g\ltimes\{\pm 1\}^g.
  \end{align*}
  
  Here the transposition $(n\quad m)$ of $S_g=\Aut(\{1,2,\dotsc,g\})$ corresponds to the element ${(n\quad m)(-n\quad {-m})}$ of $\Aut(\{\pm 1,\pm 2,\dotsc,\pm g\})$ and the element of $\{\pm 1\}^g$ which has a $-1$ in position $i$ and $1$ everywhere else corresponds to $(i \quad {-i})$.

  The affine Weyl group is $W_a=W\ltimes Y_0$ and the Iwahori-Weyl group $\tW=W\ltimes Y$ with
  \begin{align*}
\Z^{g+1}&\cong Y=\{(\nu_1,\dotsc,\nu_g,\nu_{-g},\dotsc,\nu_{-1})\in\Z^{2g} : \nu_1+\nu_{-1}=\dotsb=\nu_g+\nu_{-g}\} \\ &\supseteq Y_0=\{(\nu_1,\dotsc,\nu_g,\nu_{-g},\dotsc,\nu_{-1})\in\Z^{2g} : 0=\nu_1+\nu_{-1}=\dotsb=\nu_g+\nu_{-g}\}\cong\Z^g.
  \end{align*}

  The simple affine roots (whose walls bound the base alcove $\fa$) are
  \begin{align*}
    &1-2e_{-1}+e_0=1+2e_1-e_0, \\
    &e_{-1}-e_{-2}=e_2-e_1, e_{-2}-e_{-3}, \dotsc, e_{-(g-1)}-e_{-g}, \\
    &2e_{-g}-e_0=e_0-2e_g,
  \end{align*}
  where $e_1,\dotsc,e_g,e_{-g},\dotsc,e_{-1}\colon T\to\G_m$ are the obvious cocharacters and $e_0=e_1+e_{-1}=\dotsb=e_g+e_{-g}$.

  The reflections corresponding to the simple affine roots are
  \begin{equation*}
    ((1\quad{-1}),
    \begin{smallpmatrix}
      -1 \\ 0 \\ \vdots \\ 0 \\ 1
    \end{smallpmatrix}
    ),(-1\quad {-2})(1\quad 2),\dotsc,(-g\quad {-(g-1)})(g\quad {g-1}),(g\quad{-g}).
  \end{equation*}

  The length zero subgroup $\Omega\subseteq\tW$ is generated by $((w_0,\epsilon),y)\in (S_g\ltimes\{\pm1\}^g)\ltimes Y$, where $w_0\in S_g$ is the longest element, $\epsilon=(-1,-1,\dotsc,-1)$ and $y=(0^g,1^g)$.
\end{Remark}

\begin{Remark}\label{andereralkoven}
  One also can choose $\dotsb \subseteq p\Z_p\oplus \Z_p^{2g-1}\subseteq \Z_p^{2g}\subseteq\dotsb$ instead of $\dotsb\subseteq \Z_p^{2g-1}\oplus p\Z_p\subseteq \Z_p^{2g}\subseteq\dotsb$ as the standard lattice chain. Then the simple affine roots would be
  \begin{equation*}
    1-2e_{1}+e_0, e_{1}-e_{2}=e_2-e_1, e_{2}-e_{3}, \dotsc, e_{g-1}-e_{g}, 2e_{g}-e_0.
  \end{equation*}
\end{Remark}

\begin{Remark}\label{iw-section}
  $\tW=W\ltimes Y=N(\Q_p)/T(\Z_p)$ and $N(\Q_p) \to W\ltimes Y$ has a section $W\ltimes Y\to N(\Q_p)$, which sends $(\pi,\underline{\nu})\in W\ltimes Y$ to $T_{\underline{\nu}}P_w$, where $T_{\underline{\nu}}=
\begin{smallpmatrix}
  p^{\nu_1} & & & & \\
  & p^{\nu_2} & & &\\
  & & \ddots & & \\
  & & & p^{\nu_{-2}} & \\
  & & & & p^{\nu_{-1}}
\end{smallpmatrix}$ and $P_w$ is the permutation matrix with $P_w(e_i)=e_{w(i)}$.
\end{Remark}

\begin{Remark}\label{kottrap}
  Using the results of \cite{kottrap} we also easily can compute $\Adm(\{\mu\})$. One potential source of confusion at this point is that, due to our choice of the base alcove (cf. Remark~\ref{andereralkoven}), in our setup we need to use $\omega_i:=(0^{2g-i},1^i)$ instead of $\omega_i:=(1^i,0^{2g-i})$ (notation of \cite{kottrap}), cf. \cite[1268]{yu}. With that convention in place, we have that $x\in\tW$ is $\{\mu\}$-admissible ($\mu=(1^g,0^g)$) if and only if
  \begin{equation*}
    (0,\dotsc,0) \leq x(\omega_i) - \omega_i \leq (1,\dotsc,1) \quad\text{for all } 0\leq i <2g
  \end{equation*}
  (component-wise comparison).
\end{Remark}

\subsubsection{Lattice chains, zips, admissibility}
\label{sec:lattice-chains-zips}

\begin{Definition}\label{siegel-chain}
  Let $S$ be a $\Z_p$-scheme.
  
  A \defn{Siegel lattice chain in the weak sense on $S$ of type $J$} is a tuple $(\mathcal{V}^\bullet,\mathcal{L},\alpha_{\bullet\bullet},\theta_\bullet,\psi_\bullet)$, where
  \begin{enumerate}[(a)]
  \item for all $j\in J$, $\mathcal{V}^j$ is a vector bundle on $S$ of rank $2g$,
  \item $\mathcal{L}$ is a line bundle on $S$,
  \item for all $i,j\in J$ with $j>i$, $\alpha_{j,i}\colon \mathcal{V}^j\to \mathcal{V}^i$ is a vector bundle homomorphism, such that the $\bigl(\alpha_{j,i}\bigr)$ satisfy the obvious cocycle condition (and we also define $\alpha_{i,i}:=\id$),
  \item for all $j\in J$, $\theta_j\colon \mathcal{V}^j \xrightarrow{\sim} \mathcal{V}^{j-2g}$ is a vector bundle isomorphism such that the $\bigl(\theta_j\bigr)$ are compatible with the $\bigl(\alpha_{j,i}\bigr)$ in that $\theta_i\circ\alpha_{j,i}=\alpha_{j-2g,i-2g}\circ\theta_j$ and $\alpha_{j,j-2g}=p\cdot\theta_j$,
  \item \label{item:siegel-chain-pol} for all $j\in J$ a vector bundle isomorphism $\psi_j\colon \mathcal{V}^j \xrightarrow{\sim} (\mathcal{V}^{-j})^*\otimes\mathcal{L}$ compatible with $\bigl(\theta_j\bigr)$ and $\bigl(\alpha_{j,i}\bigr)$, such that $-\psi_j(x,y) = \psi_{-j}(y,x)$  for all $(x,y)\in \mathcal{V}^j\times \mathcal{V}^{-j}$.\footnote{By ``$(x,y)\in \mathcal{V}^j\times \mathcal{V}^{-j}$'' we of course mean that there is an open subset $U\subseteq S$ such that $(x,y)\in(\mathcal{V}^j\times \mathcal{V}^{-j})(U)$.}
  \end{enumerate}

  We also have a \defn{standard} Siegel lattice chain in the weak sense on $\Spec\Z_p$ (and hence by base change on every $\Z_p$-scheme $S$) of type $J$, namely the one given by the lattice chain $\left\{\Lambda^j \suchthat j\in J\right\}$. We can think of the standard Siegel lattice chain as either having varying $\mathcal{V}^j$ with the $\alpha_{j,i}$ being the obvious inclusion maps (e.g. (if $\{0,1\}\subseteq J$), $\mathcal{V}^1={\Z_p^{2g-1}\oplus p\Z_p} \xrightarrow{\alpha_{1,0}=\mathrm{inclusion}} \Z_p^{2g}=\mathcal{V}^0$) or as having constant $\mathcal{V}^j=\Z_p^{2g}$ with the $\alpha_{j,i}$ being diagonal matrices with all entries either $p$ or $1$ (e.g., $\mathcal{V}^1=\Z_p^{2g} \xrightarrow{\alpha_{1,0}=\diag(1,1,\dotsc,1,p)} \Z_p^{2g}=\mathcal{V}^0$). Usually the latter point of view is more convenient.

  A \defn{Siegel lattice chain on $S$ of type $J$} (or \defn{Siegel lattice chain in the strong sense on $S$ of type $J$}) then is a Siegel lattice chain in the weak sense on $S$ of type $J$ that Zariski-locally on $S$ is isomorphic to the standard chain.
\end{Definition}

\begin{Remarks}\label{symplecticity}
  \begin{enumerate}[(1)]
  \item \label{item:symplecticity-1} Let $(\mathcal{V}^\bullet,\mathcal{L},\alpha_{\bullet\bullet},\theta_\bullet,\psi_\bullet)$ be a Siegel lattice chain in the weak sense on $S$ of type $J$. Then
    $\tilde{\psi_j}:=(\tilde{\alpha}_{j,-j}^*\otimes\id_{\mathcal{L}})\circ \psi_j\colon \mathcal{V}^j\otimes\mathcal{V}^j \to \mathcal{L}$ is alternating.
    
    Here $\tilde{\alpha}_{j,-j}$ is defined as follows: Let $n\in\Z$ be maximal with $j-2gn\geq -j$. Then $\tilde{\alpha}_{j,-j}:=\alpha_{j-2gn,-j}\circ \theta_{j-2g(n-1)}\circ\dotsb\circ \theta_j$.
  \item \label{item:symplecticity-2} Note that this means that $\tilde{\psi}_j$ is (twisted) symplectic if $-j\in j+2g\Z$, i.e., if $j\in g\Z$.
  \end{enumerate}
\end{Remarks}

\begin{Proof} (of \eqref{item:symplecticity-1})
  Let $x,y\in\mathcal{V}^j$. Then
  \begin{align*}
    \tilde{\psi}_j(x,y)
    &= \psi_j(x,\tilde{\alpha}_{j,-j}(y)) \\
    &= \psi_j(x,(\alpha_{j-2gn,-j}\circ \theta_{j-2g(n-1)}\circ\dotsb\circ \theta_j)(y)) \\
    &= \psi_{2gn-j}(\alpha_{j,2gn-j}(x),(\theta_{j-2g(n-1)}\circ\dotsb\circ \theta_j)(y)) \\
    &= \psi_{2g(n-1)-j}((\theta_{2gn-j}\circ\alpha_{j,2gn-j})(x),(\theta_{j-2g(n-2)}\circ\dotsb\circ \theta_j)(y)) \\
    &= \dotsb \\
    &= \psi_{-j}((\theta_{-j+2g}\circ\dotsb\circ\theta_{2gn-j}\circ\alpha_{j,2gn-j})(x),y) \\
    &= -\psi_{j}(y,(\theta_{-j+2g}\circ\dotsb\circ\theta_{2gn-j}\circ\alpha_{j,2gn-j})(x)) \\
    &= -\tilde{\psi}_j(y,x).
  \end{align*}
\end{Proof}

\begin{Reminder}
  $\cG_K$ is the automorphism group of the standard Siegel lattice chain.
\end{Reminder}

The following definition is a generalization of \cite[Definition~3.1]{vieh-wed} in the Siegel case.

\begin{Definition}\label{def-siegel-zip}
  Let $S$ be an $\F_p$-scheme.

  A $\ocG_K$-zip over $S$ is a tuple $(\mathcal{V}^\bullet,\mathcal{L},\alpha_{\bullet\bullet},\theta_\bullet,\psi_\bullet,\mathcal{C}^\bullet,\mathcal{D}^\bullet,\varphi_{0}^{\bullet},\varphi_{1}^{\bullet},\varphi_{\mathcal{L}})$, where
  \begin{enumerate}[(a)]
  \item $(\mathcal{V}^\bullet,\mathcal{L},\alpha_{\bullet\bullet},\theta_\bullet,\psi_\bullet)$ is a Siegel lattice chain on $S$ of type $J$,
  \item for all $j\in J$, $\mathcal{C}^j\subseteq\mathcal{V}^j$ are locally direct summands of rank $g$ compatible with $\alpha_{\bullet\bullet},\theta_\bullet$, such that
    \begin{equation*}
      \mathcal{C}^j \hookrightarrow \mathcal{V}^j \overset{\psi_j}\cong (\mathcal{V}^{-j})^*\otimes\mathcal{L} \to (\mathcal{C}^{-j})^*\otimes\mathcal{L}
    \end{equation*}
    vanishes. (cf. Remark~\ref{loc-mod-siegel} for the origins of this condition.)
  \item $\mathcal{D}^\bullet\subseteq\mathcal{V}^\bullet$ satisfies the same conditions as $\mathcal{C}^\bullet\subseteq \mathcal{V}^\bullet$,
  \item $\varphi_{0}^{j}\colon (\mathcal{C}^j)^{(p)} \xrightarrow{\sim} \mathcal{V}^j/\mathcal{D}^j$ and $\varphi_{1}^{j}\colon (\mathcal{V}^j/\mathcal{C}^j)^{(p)} \xrightarrow{\sim} \mathcal{D}^j$ are isomorphisms of vector bundles compatible with $\alpha_{\bullet\bullet}$ and $\theta_\bullet$ and $\varphi_\mathcal{L}\colon \mathcal{L}^{(p)}\xrightarrow{\sim}\mathcal{L}$ is an isomorphism of line bundles, such that
    \begin{equation*}
      \begin{tikzpicture}[node distance=4cm, auto]
        \node (X) at (0,2) {$(\mathcal{C}^j)^{(p)}$};
        \node (Z) at (5,2) {$(\mathcal{V}^{-j}/\mathcal{C}^{-j})^{*,(p)}\otimes \mathcal{L}^{(p)}$};
        \node (Zs) at (5,0) {$(\mathcal{D}^{-j})^*\otimes \mathcal{L}$};
        \node (Xs) at (0,0) {$\mathcal{V}^j/\mathcal{D}^j$};
        
        \draw[->] (X) to node {$\varphi_{0}^{j}$} (Xs);
        \draw[<-] (Z) to  node {$(\varphi_{1}^{-j})^*\otimes\varphi_\mathcal{L}^{-1}$}  (Zs);
        \draw[->] (X) to node {$\psi_j^{(p)}$} (Z);
        \draw[->] (Xs) to  node{$\psi_j$}  (Zs);
      \end{tikzpicture}
    \end{equation*}
    commutes, i.e.,
    \begin{equation*}
      {\psi_j(\varphi_{0}^{j}(\_),\varphi_{1}^{-j}(\_)) = \varphi_\mathcal{L}\circ\psi_j^{(p)}(\_,\_)\colon}{(\mathcal{C}^j)^{(p)}\times (\mathcal{V}^{-j}/\mathcal{C}^{-j})^{(p)}\to \mathcal{L}^{(p)}\to \mathcal{L}}.
    \end{equation*}
  \end{enumerate}

  Since $\varphi_{\mathcal{L}}$ evidently is uniquely determined by the other data, we sometimes leave it out.

  We obtain a fibered category $\ocG_K\Zip \to \mathrm{Sch}_{\F_p}$.
\end{Definition}

\begin{Remark}\label{psi-rise}
  $\psi_j$ gives rise to isomorphisms
  \begin{align*}
    \mathcal{C}^j &\xrightarrow{\sim} (\mathcal{V}^{-j}/\mathcal{C}^{-j})^*\otimes\mathcal{L}, \\
    \mathcal{V}^j/\mathcal{C}^j &\xrightarrow{\sim} (\mathcal{C}^{-j})^*\otimes\mathcal{L}, \\
    \mathcal{D}^j &\xrightarrow{\sim} (\mathcal{V}^{-j}/\mathcal{D}^{-j})^*\otimes\mathcal{L}, \\
    \mathcal{V}^j/\mathcal{D}^j &\xrightarrow{\sim} (\mathcal{D}^{-j})^*\otimes\mathcal{L}.
  \end{align*}

  This way $\mathcal{V}^\bullet/\mathcal{C}^\bullet\oplus \mathcal{C}^\bullet$ and $\mathcal{D}^\bullet\oplus \mathcal{V}^\bullet/\mathcal{D}^\bullet$ become Siegel lattice chains in the weak(!) sense of type $J$.  The Cartier isomorphism then is an isomorphism in the category of Siegel lattice chains in the weak sense of type $J$. Over an algebraically closed field, we call the isomorphism type of the Siegel lattice chain in the weak sense $\mathcal{V}^\bullet/\mathcal{C}^\bullet\oplus \mathcal{C}^\bullet$ the \defn{Kottwitz-Rapoport type of the $\ocG_K$-zip}.
\end{Remark}

We also define a linearly rigidified version of $\ocG_K\Zip$ as follows.

\begin{Definition}
  We define the fibered category $\ocG_K\Zip^{\sim}\to \mathrm{Sch}_{\F_p}$ just like $\ocG_K\Zip$ but with the extra condition that $(\mathcal{V}^\bullet,\mathcal{L},\alpha_{\bullet\bullet},\theta_\bullet,\psi_\bullet)$ be the standard Siegel lattice chain (rather than just locally isomorphic to it).
\end{Definition}

\begin{Lemma}\label{GKZipEmb}
  We always have a closed embedding of $\ocG_K\Zip^{\sim}$ into a product of (classical) $\GL_{2g}\Zip^{\sim}$'s, and therefore $\ocG_K\Zip^{\sim}$ is a scheme.
\end{Lemma}

\begin{Proof}
  Set $J':= J\cap \{0,\dotsc,2g-1\}$. Let $\GL_{2g}\Zip^{\sim}=E_{(1^g,0^g)}\backslash (\GL_{2g}\times \GL_{2g})$ be the $\F_p$-scheme of trivialized $\GL_{2g}$-zips (so that $[\GL_{2g}\backslash\GL_{2g}\Zip^{\sim}]=\GL_{2g}\Zip$) with respect to the the cocharacter $(1^g,0^g)$, and $\prod_{j\in J'}\GL_{2g}\Zip^{\sim}$ the product of $\#J'$ copies of this scheme. On $J'$ we define $-j:=2g-j$ for $1\leq j \leq 2g-1$ and $-0:=0$.

  Then we get a monomorphism
  \begin{equation}\label{eq:GKZipEmb}
    \ocG_K\Zip^{\sim} \hookrightarrow \prod_{j\in J'}\GL_{2g}\Zip^{\sim}
  \end{equation}
  by sending $(\mathcal{C}^\bullet,\mathcal{D}^\bullet,\varphi_{0}^{\bullet},\varphi_{1}^{\bullet})$ to $\left(\mathcal{C}^j,\mathcal{D}^j,\varphi_{0}^{j},\varphi_{1}^{j}\right)_{j\in J'}$.

  The extra conditions for an element of $\prod_{j\in J'}\GL_{2g}\Zip^{\sim}$ to be in $\ocG_K\Zip^{\sim}$ are as in Definition~\ref{def-siegel-zip}:
  \begin{enumerate}[(1)]
  \item \label{item:glzip-prod-cond-1} $\mathcal{C}^\bullet,\mathcal{D}^\bullet,\varphi_{0}^{\bullet},\varphi_{1}^{\bullet}$ are compatible with the transition maps (or, to put it differently, ${(\mathcal{C}^j\oplus \mathcal{V}^j/\mathcal{C}^j)^{(p)}}\xrightarrow[\cong]{\varphi_0^j\oplus\varphi_1^j}{\mathcal{V}^j/\mathcal{D}^j\oplus \mathcal{D}^j}$ is compatible with the transition maps),
  \item \label{item:glzip-prod-cond-2} $\mathcal{C}^j \hookrightarrow \mathcal{V}^j \overset{\psi_j}\cong (\mathcal{V}^{-j})^* \to (\mathcal{C}^{-j})^*$ vanishes.
  \item \label{item:glzip-prod-cond-3} $\mathcal{D}^j \hookrightarrow \mathcal{V}^j \overset{\psi_j}\cong (\mathcal{D}^{-j})^* \to (\mathcal{D}^{-j})^*$ vanishes.
  \item \label{item:glzip-prod-cond-4} There is a (necessarily unique) isomorphism  $\varphi_\mathcal{L}\colon \mathcal{L}^{(p)}\xrightarrow{\sim}\mathcal{L}=\mathcal{O}_S$ of line bundles, such that
    \begin{equation*}
      \begin{tikzpicture}[node distance=4cm, auto]
        \node (X) at (0,2) {$(\mathcal{C}^j)^{(p)}$};
        \node (Z) at (5,2) {$(\mathcal{V}^{-j}/\mathcal{C}^{-j})^{*,(p)}$};
        \node (Zs) at (5,0) {$(\mathcal{D}^{-j})^*$};
        \node (Xs) at (0,0) {$\mathcal{V}^j/\mathcal{D}^j$};
        
        \draw[->] (X) to node {$\varphi_{0}^{j}$} (Xs);
        \draw[<-] (Z) to  node {$(\varphi_{1}^{-j})^*\otimes\varphi_\mathcal{L}^{-1}$}  (Zs);
        \draw[->] (X) to node {$\psi_j^{(p)}$} (Z);
        \draw[->] (Xs) to  node{$\psi_j$}  (Zs);
      \end{tikzpicture}
    \end{equation*}
    commutes.
  \end{enumerate}

  We claim that the conditions are closed on $\prod_{j\in J'}\GL_{2g}\Zip^{\sim}$  (hence the monomorphism is a closed immersion).

  To see this, we recall the construction of the scheme $\GL_{2g}\Zip^{\sim}$ as executed in \cite[\nopp (3.10), (3.11), (4.3)]{mw}.

  Recall our notational convention regarding the parabolic subgroup associated with a cocharacter $\chi$ from Definition~\ref{pappzhu}. As in \cite{mw}, we denote by $\mathrm{Par}_\chi$ the scheme of parabolic subgroups of type $\chi$.

  There is a group scheme $H$ defined by the cartesian diagram
  \begin{equation*}
    \begin{tikzpicture}[node distance=4cm, auto]
      \node (LO) at (0,3) {$H$};
      \node (RO) at (5,3) {$\mathcal{P}_{((-1)^g,0^g)}/\mathcal{U}_{((-1)^g,0^g)}$};
      \node (sq) at (2.5,1.5) {$\square$};
      \node (LU) at (0,0) {$\mathrm{Par}_{((-1)^g,0^g)}\times\mathrm{Par}_{(1^g,0^g)}$};
      \node (RU) at (5,0) {$\mathrm{Par}_{((-1)^g,0^g)}$};
      
      \draw[->] (LO) to  (LU);
      \draw[->] (RO) to  (RU);
      \draw[->] (LO) to  (RO);
      \draw[->] (LU) to node {$(\;)^{(p)}\circ \pr_1$} (RU);
    \end{tikzpicture}
  \end{equation*}
  where $\mathcal{P}_{((-1)^g,0^g)}\to \mathrm{Par}_{((-1)^g,0^g)}$ is the universal parabolic group scheme and $\mathcal{U}_{((-1)^g,0^g)}$ its unipotent radical, such that $\GL_{2g}\Zip^{\sim}$ is an $H$-Zariski torsor over $\mathrm{Par}_{((-1)^g,0^g)}\times\mathrm{Par}_{(1^g,0^g)}$, where $\GL_{2g}\Zip^{\sim}\to \mathrm{Par}_{((-1)^g,0^g)}\times\mathrm{Par}_{(1^g,0^g)}$ is given by $(C,D,\varphi_0,\varphi_1)\mapsto (C,D)$.

  Clearly, compatibility of $\mathcal{C}^\bullet,\mathcal{D}^\bullet$ with the transition maps is a closed condition on $\prod_{j\in J'}\mathrm{Par}_{((-1)^g,0^g)}\times\mathrm{Par}_{(1^g,0^g)}$ and then also on  $\prod_{j\in J'}\GL_{2g}\Zip^{\sim}$. Similar for the conditions \eqref{item:glzip-prod-cond-2} and \eqref{item:glzip-prod-cond-3}.

  Locally, we can choose complements (not necessarily compatible with the transition maps) and then $\varphi_\bullet^j$ yield sections $g^j$ of $\GL_{2g}$ as in \cite[\nopp definition of $g\in G(S)$ in the proof of (4.3)]{mw}. The $g^j$ are well-defined up to $\mathcal{U}_{((-1)^g,0^g)}^{(p)}\times\mathcal{U}_{(1^g,0^g)}$, and we want them to be compatible with the transition maps coming from the Siegel lattice chains in the weak sense $\mathcal{C}^j\oplus\mathcal{V}^j/\mathcal{C}^j$ and $\mathcal{V}^j/\mathcal{D}^j\oplus\mathcal{D}^j$, respectively. With our complements in place, these transition maps correspond to maps $\mathcal{V}^j\to\mathcal{V}^{j-n}$. The question of whether $g^j$ is compatible with these maps is independent of the choice of complements (basically because the transition maps $\mathcal{V}^j\to\mathcal{V}^{j-n}$ depend on the choice of complements similar to how $g^j$ depends on that choice).

  So in effect we can view the conditions on $\varphi_0^\bullet,\varphi_1^\bullet$ of \eqref{item:glzip-prod-cond-1} as closed conditions on $\prod_{j\in J'}\GL_{2g}\Zip^{\sim\sim}$, where $\GL_{2g}\Zip^{\sim\sim}\to\GL_{2g}\Zip^{\sim}$ (an fpqc quotient map) additionally comes with complementary spaces of $C,D$ ($\GL_{2g}\Zip^{\sim\sim}=\tilde{X}_\tau$ in the notation of \cite[\nopp proof of (4.3)]{mw}).

  We also can reformulate condition \eqref{item:glzip-prod-cond-4} in those terms.
\end{Proof}

\begin{Corollary}
  $\ocG_K\Zip$ is the algebraic quotient stack $[\ocG_K\backslash\ocG_K\Zip^{\sim}]$.

  Here by definition an element $\phi\in\ocG_K$ acts on $\ocG_K\Zip^{\sim}$ by replacing $(\mathcal{C}^\bullet,\mathcal{D}^\bullet,\varphi_{0}^{\bullet},\varphi_{1}^{\bullet},\varphi_{\mathcal{L}})$ by \linebreak[3] $(\phi\mathcal{C}^\bullet,\phi\mathcal{D}^\bullet,\phi\varphi_{0}^{\bullet}\phi^{-(p)},\phi\varphi_{1}^{\bullet}\phi^{-(p)},\varphi_{\mathcal{L}})$.
\end{Corollary}

\begin{Definition}
  We let an element $(X,Y)\in\ocG_K\times\ocG_K$ act on $\ocG_K\Zip$ by replacing $(\mathcal{V}^\bullet,\mathcal{L},\alpha_{\bullet\bullet},\theta_\bullet,\psi_\bullet,\mathcal{C}^\bullet,\mathcal{D}^\bullet,\varphi_{0}^{\bullet},\varphi_{1}^{\bullet})$ by
  \begin{equation*}
    (\mathcal{V}^\bullet,\mathcal{L},\alpha_{\bullet\bullet},\theta_\bullet,\psi_\bullet,
    X\mathcal{C}^\bullet,Y\mathcal{D}^\bullet,
    Y\varphi_{0}^{\bullet}X^{-(p)},Y\varphi_{1}^{\bullet}X^{-(p)}).
  \end{equation*}
\end{Definition}

\begin{Notation}
  Let $\sS_K\to\Spec\Z_p$ be the integral model of the Siegel Shimura variety of level $K$ (where $K=K_pK^p$ with $K^p$ sufficiently small), and recall $\tsS_K$ from Section~\ref{sec:siegel-case-loc-mod}. Moreover, define $\osS_K:=\sS_K\otimes_{\Z_p}\F_p$ and $\tosS_K:=\tsS_K\otimes_{\Z_p}\F_p$. So $\tosS_K\to\osS_K$ is a $\ocG_K$-torsor.
\end{Notation}

\begin{Remark}
  We have morphisms $\tosS_K \to \ocG_K\Zip^{\sim}$ (take first de Rham cohomology with Frobenius and Verschiebung) and $\ocG_K\Zip^{\sim} \to \overline{M}^\mathrm{loc}_K$ (take the $\mathcal{C}^\bullet$-filtration) and therefore
  \begin{equation*}
    \osS_K \to \ocG_K\Zip \to [\ocG_K\backslash \overline{M}^\mathrm{loc}_K].
  \end{equation*}
\end{Remark}

\begin{Remark}
  In particular, $\ocG_K\Zip$ has a Kottwitz-Rapoport stratification, which agrees with the notion of Kottwitz-Rapoport type as defined in Remark~\ref{psi-rise}.

  For $w\in \KR(K,\{\mu\})$ denote the associated Kottwitz-Rapoport stratum by $\ocG_K\Zip_w$, i.e.,
  we interpret $w$ as a $\bar\F_p$-valued point of $[\ocG_K\backslash \overline{M}^\mathrm{loc}_K]$ and form $\ocG_K\Zip_w$ as a fiber product.
\end{Remark}

\begin{Construction}\label{std-zeug}
  Fix $w\in \Adm(\{\mu\})^K\subseteq\tW$ (so that $W_KwW_K\in\KR(K,\{\mu\})$). We define a standard $\ocG_K$-zip of KR type $W_KwW_K$.

  Using Remark~\ref{iw-section}, we interpret $w$ as an element of $N(\Q_p)\subseteq G(\Q_p)$. The admissibility condition implies that we can interpret it as an endomorphism $w^\bullet$ of the standard lattice chain $\mathcal{V}^\bullet$ over $\Z_p$.\footnote{Take up the second point of view described in Definition~\ref{siegel-chain} regarding $\mathcal{V}^\bullet$. Define $\underline{\nu}^{(0)}:=\underline{\nu}$,
  $\underline{\nu}^{(1)}:=\underline{\nu}+
  \begin{smallpmatrix}
    0 \\ \vdots \\ 0 \\ 0 \\ -1
  \end{smallpmatrix}+w
  \begin{smallpmatrix}
    0 \\ \vdots \\ 0 \\ 0 \\ 1
  \end{smallpmatrix}$,
  $\underline{\nu}^{(2)}:=\underline{\nu}+
  \begin{smallpmatrix}
    0 \\ \vdots \\ 0 \\ -1 \\ -1
  \end{smallpmatrix}+w
  \begin{smallpmatrix}
    0 \\ \vdots \\ 0 \\ 1 \\ 1
  \end{smallpmatrix}$, and so on. Then $w^j=T_{\underline{\nu}^{(j)}}P_w$ for $0\leq j<2g$.

  From the formulation of the admissibility condition as in Remark~\ref{kottrap}, we see that $w\in\Adm(\{\mu\})^K$ is equivalent to the condition that $\underline{\nu}^{(j)}$ be a permutation of $(1^g,0^g)$ for all relevant $j$.}

We denote the standard Siegel lattice chain over $\Z_p$ by $\sV^\bullet$ and its base change to $\F_p$ by $\mathcal{V}^\bullet$. Define $\sC_w^\bullet:=pw^{\bullet,-1}\sV^\bullet$ and $\sD_w^\bullet:=\sigma(w^\bullet)\sV^\bullet$. Then $\mathcal{C}_w^\bullet:=\sC_w^\bullet\otimes\F_p=\ker(w^\bullet\colon \mathcal{V}^\bullet\to \mathcal{V}^\bullet)$, so $(\mathcal{V}^\bullet/\mathcal{C}_w^\bullet)^{(p)}\xrightarrow{\sim} \mathcal{D}_w^\bullet:=\sD_w^\bullet\otimes\F_p$ via $\sigma(w^\bullet)$ and $(\mathcal{C}_w^\bullet)^{(p)}\xrightarrow{\sim} \mathcal{V}^\bullet/\mathcal{D}_w^\bullet$ via $p^{-1}\sigma(w^\bullet)$.

  This defines a standard element $\tStd(w)$ of $\ocG_K\Zip^{\sim}(\F_p)$ and a standard element $\Std(w)$ of $\ocG_K\Zip_w(\F_p)$.
\end{Construction}

\begin{DefinitionRemark} \textnormal{(See also \cite[Lemma~3.3.2]{SYZnew}.)}
  $\cG_w:=\Aut(\sC_w^\bullet\subseteq \sV^\bullet)$ is a Bruhat-Tits group scheme with generic fiber $G_{\Q_p}$ and $\bZ_p$-points $\bK\cap w^{-1}\bK w$; and similarly for $\cG_{\sigma(w)^{-1}}:=\Aut(\sD_w^\bullet\subseteq \sV^\bullet)$ with $\bK\cap \sigma(w)\bK\sigma(w)^{-1}$.
\end{DefinitionRemark}

\begin{Definition}
  We keep $w\in \Adm(\{\mu\})^K\subseteq\tW$ fixed and define $\tE_w\subseteq\ocG_K\times\ocG_K$ to be the stabilizer of $\tStd(w)$.

  So $\tE_w$ consists of those $(X^\bullet,Y^\bullet)\in\ocG_K\times\ocG_K$ such that $X^\bullet\mathcal{C}_w^\bullet=\mathcal{C}_w^\bullet$, $Y^\bullet\mathcal{D}_w^\bullet=\mathcal{D}_w^\bullet$, and $Y^\bullet\circ\varphi_j^\bullet\circ X^{\bullet,-(p)}=\varphi_j^\bullet$ for $j=0,1$.

  In the notation of \cite[Lemma~3.3.2]{SYZnew} we have
  \begin{equation}\label{eq:verallg-zip-gp}
    \tE_w=\ocG_w\times_{\ocG_w^{L,(p)}}\ocG_{\sigma(w)^{-1}}.
  \end{equation}

  Here $\ocG_w^{L}$ is the image of $\ocG_w$ in $\mathrm{DiagAut}(\mathcal{C}_w^\bullet\oplus \mathcal{V}^\bullet/\mathcal{C}_w^\bullet)$ (the automorphisms of $\mathcal{C}_w^\bullet\oplus \mathcal{V}^\bullet/\mathcal{C}_w^\bullet$ respecting both $\mathcal{C}_w^\bullet$ and $\mathcal{V}^\bullet/\mathcal{C}_w^\bullet$).

  The orbit of $\tStd(w)$ in $\ocG_K\Zip^{\sim}$ is the fppf quotient $(\ocG_K\times\ocG_K)/\tE_w$, cf. \cite[\nopp II, \S\,5, no.~3]{DGen}.
\end{Definition}

\begin{Lemma}\label{diagsmorris}
  We have commutative diagrams
  \begin{equation*}
    \begin{tikzpicture}[node distance=4cm, auto]
      \node (X) at (0,2) {$\ocG_w$};
      \node (Z) at (5,2) {$\ocG$};
      \node (Zs) at (5,0) {$\ocG^\mathrm{rdt}$};
      \node (Xs) at (0,0) {$\bar{P}_{J_1}$};
      
      \draw[->>] (X) to (Xs);
      \draw[->>] (Z) to  (Zs);
      \draw[right hook->] (X) to  (Z);
      \draw[right hook->] (Xs) to  (Zs);
    \end{tikzpicture}
  \end{equation*}
  and
  \begin{equation*}
    \begin{tikzpicture}[node distance=4cm, auto]
      \node (X) at (0,2) {$\ocG_{\sigma(w)^{-1}}$};
      \node (Z) at (5,2) {$\ocG$};
      \node (Zs) at (5,0) {$\ocG^\mathrm{rdt}$};
      \node (Xs) at (0,0) {$\bar{P}_{\sigma'(J_1)}$};
      
      \draw[->>] (X) to (Xs);
      \draw[->>] (Z) to  (Zs);
      \draw[right hook->] (X) to  (Z);
      \draw[right hook->] (Xs) to  (Zs);
    \end{tikzpicture}
  \end{equation*}
  and
  \begin{equation*}
    \begin{tikzpicture}[node distance=4cm, auto]
      \node (X) at (0,2) {$\ocG_w^L$};
      \node (Z) at (5,2) {$\ocG$};
      \node (Zs) at (5,0) {$\ocG^\mathrm{rdt}$};
      \node (Xs) at (0,0) {$\bar{L}_{J_1}$.};
      
      \draw[->>] (X) to (Xs);
      \draw[->>] (Z) to  (Zs);
      \draw[right hook->] (X) to  (Z);
      \draw[right hook->] (Xs) to  (Zs);
    \end{tikzpicture}
  \end{equation*}
\end{Lemma}

\begin{Proof}
  This follows from Proposition~\ref{propmorris}.
\end{Proof}

\begin{Lemma}\label{Ebild}
  The image of $\tE_w$ under $\ocG\times\ocG \to \ocG^\mathrm{rdt}\times\ocG^\mathrm{rdt}$ is $E_{\mathcal{Z}_w}$.
\end{Lemma}

\begin{Proof}
  This follows from Lemma~\ref{diagsmorris}.
\end{Proof}

\begin{Lemma} \label{orbits-kr} Assume $0\in J$.
  The $\ocG_K\times\ocG_K$-orbit of $\tStd(w)$ for $w\in\Adm(\{\mu\})^K$ depends only on $W_KwW_K$.
\end{Lemma}

\begin{Proof}
  Let $x,y\in W_K \subseteq W$. As above we get endomorphisms $x^\bullet,y^\bullet$ of $\mathcal{V}^\bullet$, which in this case are in fact automorphisms. Now $\tStd(w)=((y^\bullet)^{-1},\sigma(x^\bullet))\cdot\tStd(w)$.
\end{Proof}

\begin{Definition}
  Define $\ocG_K\AdmZip^{\sim}$ to be the union of the $\ocG_K\times\ocG_K$-orbits of the standard zips $\tStd(w)$ for $w\in\Adm(\{\mu\})^K$. Here an orbit by definition is the image of the orbit map endowed with the reduced subscheme structure, and---as we prove just below---the union of orbits just referred to is a closed subset, which we again endow with the reduced subscheme structure.

  Define $\ocG_K\AdmZip:=[\ocG_K\backslash\ocG_K\AdmZip^{\sim}]\subseteq [\ocG_K\backslash\ocG_K\Zip^{\sim}]=\ocG_K\Zip$.
\end{Definition}

\begin{Lemma}\label{admzip-closed}
  $\ocG_K\AdmZip^{\sim}$ is a closed subset of $\ocG_K\Zip^{\sim}$.
\end{Lemma}

\begin{Proof}
  This being a purely topological question, we may freely pass to perfections, which will be convenient since DieudonnÃ© theory is simpler over perfect rings. By ``perfection'' we mean the inverse perfection in the terminology of \cite[Section~5]{perfection}.

  Consider therefore $(\ocG_K\Zip^{\sim})^\mathrm{perf}$ as a sheaf on $\mathrm{Perf}_{\F_p}$, the fpqc site of affine perfect $\F_p$-schemes. Again denoting the standard Siegel lattice chain over $\Z_p$ by $\sV^\bullet$ and its base change to $\F_p$ by $\mathcal{V}^\bullet$, we can describe the elements of $(\ocG_K\Zip^{\sim})^\mathrm{perf}(R)=\ocG_K\Zip^{\sim}(R)$, where $R$ is a perfect $\F_p$-algebra as being given by
  \begin{quote}
    homomorphisms $\mathcal{V}_R^{\bullet,(p)} \xrightarrow{F^\bullet} \mathcal{V}_R^\bullet \xrightarrow{V^\bullet} \mathcal{V}_R^{\bullet,(p)}$ such that $\ker(F^\bullet)=:\mathcal{C}^{\bullet,(p)}=\im(V^\bullet)$ and $\im(F^\bullet)=:\mathcal{D}^\bullet=\ker(V^\bullet)$ and $\psi_j(F^j\_,\_)=u\sigma(\psi_j(\_,V^{-j}\_))$ for some $u\in R^\times$ and $\mathcal{C}^{\bullet,(p)},\mathcal{D}^\bullet$ have the same rank (namely $g$).
  \end{quote}
  To see that $\mathcal{C}^{\bullet,(p)},\mathcal{D}^\bullet$ are direct summands of $\mathcal{V}_R^{\bullet,(p)},\mathcal{V}_R^\bullet$ (which makes the last part of the characterization given above meaningful), one argues as in \cite[Lemma~2.4]{laurel} (since both are finitely presented, it is enough to show flatness and to that end, one looks at the fiber dimensions).

  Define a presheaf $\mathcal{X}$ on $\mathrm{Sch}_{\Z_p}$ in the same way but for the following changes: $\mathcal{V}^\bullet$ is replaced by $\sV^\bullet$, and we impose the condition that both compositions $F^\bullet\circ V^\bullet$ and $V^\bullet\circ F^\bullet$ are multiplication by $p$, and the $\ker=\im$-conditions are only required to hold modulo $p$. We also slightly reformulate these $\ker=\im$-conditions: We impose the condition that the  reductions $\bar F^\bullet,\bar V^\bullet$ be fiberwise of rank $g$ over $R/p$. (Note that the argument that $\mathcal{C}^{\bullet,(p)},\mathcal{D}^\bullet$ are direct summands only works over reduced rings.)

  Then $\mathcal{X}$ is a separated $\Z_p$-scheme. To see this, we build it up from scratch as follows.
  $\End(\sV^j)$ obviously is a $\Z_p$-scheme (an affine space), hence so is $\Hom(\sV^{j,(p)},\sV^j)$ since
  $\sV_j^{(p)}\cong \sV_j$.
  $\Hom(\sV^{\bullet,(p)},\sV^\bullet)$ is a locally closed subscheme of a finite product of such schemes.
  Homomorphisms $\sV^{\bullet,(p)} \xrightarrow{F^\bullet} \sV^\bullet \xrightarrow{V^\bullet} \sV^{\bullet,(p)}$ such that both compositions are multiplication by $p$ form a closed subscheme $\mathcal{X}'$ of $\Hom(\sV^{\bullet,(p)},\sV^\bullet) \times \Hom(\sV^\bullet,\sV^{\bullet,(p)})$.
  In the special fiber $\mathcal{X}'_{\F_p}$ we now consider the $\ker=\im$-conditions and show that they define an open subscheme $\bar{\mathcal{X}}''$. Then $\mathcal{X}=\mathcal{X}'\times_{\mathcal{X}'_{\F_p}}\bar{\mathcal{X}}''$.
  Indeed, the extra conditions are that all $F^\bullet,V^\bullet$ have some non-vanishing $g$-minor---evidently open conditions. 

  The upshot is that we defined a  $\Z_p$-scheme $\mathcal{X}$ such that $(\mathcal{X}\times_{\Z_p}\F_p)^\mathrm{perf}=(\ocG_K\Zip^{\sim})^\mathrm{perf}$ and such that we have an obvious morphism $\tsS_K\to \mathcal{X}$, which takes a principally polarized isogeny chain of abelian schemes to the evaluation of the DieudonnÃ© crystal on the trivial thickening.\footnote{This makes use of the crystalline-de Rham comparison to make a trivialization of the de Rham cohomology into a trivialization of the crystalline cohomology.}

  Observe that $\mathcal{X}$ also has a natural $\mathcal{G}_K\times\mathcal{G}_K$-action: We interpret $\mathcal{G}_K$ as $\Aut(\sV^\bullet)$ and the action of $(X^\bullet,Y^\bullet)$ transforms $(F^\bullet,V^\bullet)$ into $(Y^\bullet\circ F^\bullet \circ X^{\bullet,-(p)},  X^{(p)}\circ V^\bullet\circ Y^{\bullet,-1})$. The identity $(\mathcal{X}\times_{\Z_p}\F_p)^\mathrm{perf}=(\ocG_K\Zip^{\sim})^\mathrm{perf}$ is an identity of $\ocG_K^\mathrm{perf}\times\ocG_K^\mathrm{perf}$-varieties.

  Now we claim that $\ocG_K\AdmZip^{\sim} = \left(\mathcal{X}_{\F_p}\times_{\mathcal{X}}\overline{\mathcal{X}_{\Q_p}}\right)^\mathrm{perf}$ topologically, where $\overline{\mathcal{X}_{\Q_p}}$ is the flat closure of the generic fiber in $\mathcal{X}$. This of course implies
  $\ocG_K\AdmZip^{\sim}\subseteq \ocG_K\Zip^{\sim}$ being closed.

  Both sets are constructible, so it suffices to check it on a very dense subset, say the $\bar\F_p$-valued points.

  Using Lemmas~\ref{closed-generizations}~and~\ref{closed-generizations-realization}, we see that $(\mathcal{X}_{\F_p}\times_{\mathcal{X}}\overline{\mathcal{X}_{\Q_p}})(\bar\F_p)$ consists precisely of those elements $\bar{x}\in\ocG_K\Zip^{\sim}(\bar\F_p)$ such that there exists a finite field extension $L/\bQ_p$ and a point $x\in\mathcal{X}(\mathcal{O}_L)$ lifting $\bar x$. (We'll also say that $\bar x$ is \defn{liftable} in this situation.)

  Since $\cG_K$ is flat over $\Z_p$, this liftability condition for $\mathcal{G}_K$ (in lieu of $\mathcal{X}$) is always satisfied. Consequently, $(\mathcal{X}_{\F_p}\times_{\mathcal{X}}\overline{\mathcal{X}_{\Q_p}})(\bar\F_p)$ is stable under the $\ocG_K\times\ocG_K$-action.

  Also, the standard zips clearly are liftable. Thus, $(\mathcal{X}_{\F_p}\times_{\mathcal{X}}\overline{\mathcal{X}_{\Q_p}})(\bar\F_p)\supseteq\ocG_K\AdmZip^{\sim}(\bar\F_p)$.

  For the converse inclusion, there are injective maps from $\mathcal{X}(\mathcal{O}_L)$ to $\mathcal{X}(L)$ to $\cG_K(L)$ such that the corresponding Schubert cell (in the local model) is indexed by the image mod $\cG_K(\mathcal{O}_L)\times \cG_K(\mathcal{O}_L)^\mathrm{op}$, cf. Proposition~\ref{bt-decomp}.\footnote{Note that $\cG_K(\mathcal{O}_L)\backslash \cG_K(L) /\cG_K(\mathcal{O}_L)\cong W_K\backslash\tW/W_K$ for every strictly henselian discretely valued field $L$ by \cite[Prop.~8]{hainesrap}. (And also in the construction of $\tW$ and $W_K$ any such field, not just $L=\bQ_p$, can be used.)} This proves it since we know which Schubert cells belong to the local model.
\end{Proof}

\begin{Remark}
  Regarding the orbit closure relations for $\ocG_K\AdmZip^{\sim}$, let us point out that $\ocG_K\AdmZip^{\sim}\to \bar{M}^\mathrm{loc}_K$ is $\ocG_K\times\ocG_K$-equivariant, where the action of $\ocG_K\times\ocG_K$ on $M^\mathrm{loc}_K$ factors through the first projection map, and this map is a bijection on orbits. Writing $w'\preceq w$ if $(\ocG_K\times\ocG_K)\cdot \Std(w')\subseteq \overline{(\ocG_K\times\ocG_K)\cdot \Std(w)}$, it follows from these observations that $w'\leq w$ implies $w'\preceq w$.  Here $\leq$ is the Bruhat order on $W_K\backslash \tW/W_K$ as explained in \cite[section~4.2]{prs}.
\end{Remark}

It appears reasonable to suspect that $\preceq$ and $\leq$ in fact agree.

\begin{Conjecture}
  The closure of ${(\ocG_K\times\ocG_K)\cdot\tStd(w)}$ is given by the disjoint union of ${(\ocG_K\times\ocG_K)\cdot\tStd(w')}$ for $w'\leq w$.
\end{Conjecture}

\begin{Lemma}\label{factor-thru-admzip}
  The map $\osS_K\to \ocG_K\Zip$ factors through $\ocG_K\AdmZip$.
\end{Lemma}

\begin{Proof}
  It is sufficient to check this on $k=\bar\F_p$-valued points.

  The map $\osS_K(k)\to \ocG_K\Zip(k)$ factors through $\Upsilon_K\colon \osS_K(k) \to \bigcup_{w\in \KR(K,\{\mu\})} \bK w \bK/\bK_\sigma$ with $\bK w \bK/\bK_\sigma\to\ocG_K\Zip(k)$ given by sending $xwy$ to $(\bar{y}^{-1},\sigma(\bar{x}))\cdot\Std(w)$ (similar to Lemma~\ref{orbits-kr}).
\end{Proof}

\subsubsection{An explicit description of \texorpdfstring{$\ocG_K^\mathrm{rdt}$}{GK\textasciicircum{}rdt}}
\label{sec:GKrdt-siegel}

In order to get a better feeling for the passage from $\ocG_K$ to the maximal reductive quotient $\ocG_K^\mathrm{rdt}=\ocG_K/R_u\ocG_K$ (with $R_u\ocG_K$ being the unipotent radical of $\ocG_K$), which is key in the definition of the EKOR stratification, we describe $\ocG_K^\mathrm{rdt}$ in explicit, linear-algebraic terms in the Siegel case.

Let $(\mathcal{V}^\bullet,\mathcal{L},\alpha_{\bullet\bullet},\theta_\bullet,\psi_\bullet)$ be the standard Siegel lattice chain on $S$ of type $J$. Assume $0\in J$. In what follows, we sometimes use $j$ as a shorthand for $\mathcal{V}^j$.

By a \defn{symmetric transition map}, we mean a transition map from $j'$ to $j''$, where $n\in\Z$, $j',j''\in J$, $ng\geq j'\geq j''> (n-2)g$, and $j' + j''\in 2g\Z$. We will also call this the symmetric transition map of $(j',n)$ (or of $j'$ if $n$ doesn't matter).
  
By a \defn{one-sided transition map}, we mean a transition map from $j'$ to $j''$, where $n\in\Z$, $j',j''\in J$, $ng\geq j'\geq j''\geq (n-1)g$. Call it right-anchored if $j'=ng$ and left-anchored if $j''=(n-1)g$. We then also speak of the right-anchored transition map of $j''$ and the left-anchored transition map of $j'$, respectively. 

The kernels of the symmetric transition maps are symplectic subbundles of $\mathcal{O}_S^{2g}$ (even of the form $\mathcal{O}_S^I$, where $I\subseteq\{\pm 1,\dotsc,\pm g\}$ is symmetric (i.e., $-I=I$)), and the kernels of the one-sided transition maps are totally isotropic subbundles (even of the form $\mathcal{O}_S^I$, where $I\subseteq \{1,\dotsc,g\}$ or $I\subseteq\{-1,\dotsc,-g\}$).

Let $\mathcal{O}_S^{I_j}$ be the kernel of the symmetric transition map of $j$. Then $I_j \sqcup I_{-j}=\{\pm 1,\dotsc,\pm g\}$.

Every kernel of a one-sided transition map is a subbundle of a kernel of an anchored transition map inside of which it is complemented by the kernel of another one-sided transition map.

The kernel of the left-anchored transition map of $j$ is a subbundle of the kernel of the symmetric transition map of $-j$ inside of which it is complemented by the kernel of the right-anchored transition map of $-j$. Likewise, the kernel of the right-anchored transition map of $j$ is a subbundle of the kernel of the symmetric transition map of $j$ inside of which it is complemented by the kernel of the left-anchored transition map of $-j$.

Now consider the standard symplectic bundle $\mathcal{O}_S^{2g}$ together with the kernels of all the symmetric transition maps and all the one-sided transition maps. So we have a symplectic bundle with a bunch of symplectic subbundles coming in complementary pairs, some of which come with a further decomposition into complementary Lagrangians, some of which come with further decompositions into complementary subbundles (of course still totally isotropic). We will also call these kernels \defn{distinguished subspaces}. 

Below we prove that $\ocG_K^\mathrm{rdt}$ is the automorphism group scheme $\mathcal{A}$ of these data. Clearly, $\mathcal{A}$ is reductive; in fact it is a Levi subgroup of a parabolic of $\GSp_{2g}$.

We have a map $\ocG_K\to\mathcal{A}$; the image of an $S$-point $f^\bullet$ under $\ocG_K\to\mathcal{A}$ on the kernel of a transition map starting at $j$ is given by $f^j$. Note that $f^j=\tau\circ f^j$ on $\ker(\tau)$ for every transition map $\tau$ starting at $j$.

$\ocG_K\to\mathcal{A}$ has a natural section $\mathcal{A}\to\ocG_K$, where in the image all the $f^j$ are the same as automorphisms of $\mathcal{O}_S^{2g}$. (This is well-defined!)

\begin{Proposition}\label{GKrdt}
  $\mathcal{A}=\ocG_K^\mathrm{rdt}$.
\end{Proposition}

\begin{Proof}
  Let us show that $\mathcal{K}:=\ker(\ocG_K\to\mathcal{A})$ is unipotent. Consider $\ocG_K$ as a subgroup of $\prod_{j\in J/2g\Z}\GL_{2g}\subseteq \GL_{N}$. We claim that said kernel is contained in $\prod_{j\in J/2g\Z} U^{(j)}$, $U^{(j)}$ being a conjugate of the standard unipotent subgroup $
  \begin{smallpmatrix}
    1 & \ast & \ast & \dotsb & \ast \\
    & 1 & \ast & \dotsb & \ast \\
     &  & \ddots & \dotsb & \vdots
  \end{smallpmatrix}
$ of $\GL_{2g}$.

  Indeed, say $f^\bullet$ is in the kernel. Then $f^j$ acts as the identity on the kernel of the symmetric transition map of $j$ and $f^{-j}$ acts as the identity on the kernel of the symmetric transition map of $-j$. On the image of the symmetric transition map $\tau_j$ of $j$, $f^{-j}$ agrees with $\tau_j\circ f^j$. Note that $\im(\tau_j)=\ker(\tau_{-j})$. So $\tau_j\circ f^j$ is the identity on $\ker(\tau_{-j})$. Hence, if $x\in \ker(\tau_{-j})$, then $x=\tau_j(x)$ and $f^j(x)\equiv x \mod \ker(\tau_j)$. Thus with respect to the decomposition $\ker(\tau_j)\oplus\ker(\tau_{-j})$, $f^j$ is of the form $
  \begin{pmatrix}
    1 & \ast \\
    & 1
  \end{pmatrix}$.

  Now we have $\ocG_K=\mathcal{A}\ltimes\mathcal{K}$, in particular $\ocG_K\cong\mathcal{A}\times_{\F_p} \mathcal{K}$ as schemes. Since both $\ocG_K$ and $\mathcal{A}$ are reduced and connected, so is $\mathcal{K}$.

  All in all, we see that  $\mathcal{A}$ is indeed $\ocG_K^\mathrm{rdt}$ and $\mathcal{K}=R_u\ocG_K$ is the unipotent radical of $\ocG_K$.
\end{Proof}

\begin{Example}\label{GKrdt-ex}
  \begin{itemize}
  \item If $J=\Z$, then $\ocG_K^\mathrm{rdt}=\G_m^{g+1}$ is the standard maximal torus of $\GSp_{2g}$.
  \item If $g=2$ and $J=2\Z$, then $\ocG_K^\mathrm{rdt}$ is the automorphism group of the standard twisted symplectic space $\F_p^4$ with its standard Lagrangian decomposition, i.e., $\ocG_K^\mathrm{rdt}\cong\GL_2\times \G_m$.
  \item If $g=2$ and $J/2g\Z=\{-1,0,1\}$, then $\ocG_K^\mathrm{rdt}$ is the automorphism group of the standard twisted symplectic space $\F_p^4$ with its standard decomposition in twisted symplectic subspaces and the totally isotropic  rank-$1$ subbundles generated by $e_{\pm 1}$, i.e., $\ocG_K^\mathrm{rdt}\cong\GL_2\times\G_m$.
  \item Let $g=8$. We have the local Dynkin diagram
    \begin{center}
      \includegraphics[height=32pt
      ]{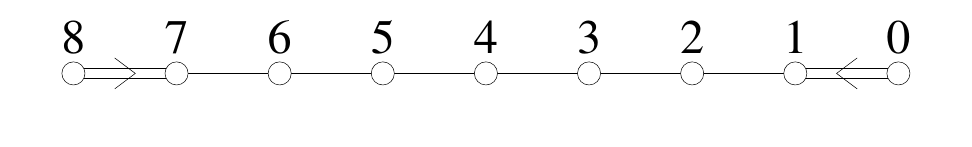}
    \end{center}
    where we labelled the simple affine roots as follows: $1-2e_{-1}+e_0$ is labelled 0, $e_{-i}-e_{-(i+1)}$ is labelled $i$ for $1\leq i\leq 7$, and $2e_{-8}-e_0$ is labelled 8.

    Consider $J/2g\Z=\{0,\pm 3,\pm 5\}$. Then the Dynkin diagram of $\ocG_K^\mathrm{rdt}$ should (according to \cite[\nopp 3.5.1]{tits}) be the one we get by removing $0,3,5$ and the adjacent edges. So we expect something along the lines\footnote{i.e., having the same Dynkin diagram as} of $\GSp(6)\times\GL(2)\times\GL(3)$.

    We have the following (bases of) kernels of symmetric transition maps:
    \begin{equation*}
      \{\pm 1, \pm 2, \pm 3\},\{\pm 4, \pm 5, \pm 6, \pm 7, \pm 8\},\{\pm 1, \pm 2, \pm 3, \pm 4, \pm 5\},\{\pm 6, \pm 7, \pm 8\},
    \end{equation*}
    and the following kernels of one-sided transition maps:
    \begin{align*}
      \{-3, -2, -1\},
      \{-5, -4\},
      \{-5, -4, -3, -2, -1\}, \\
      \{4, 5\},
      \{1, 2, 3, 4, 5\},
      \{1, 2, 3\}.
    \end{align*}

    So an element $A$ of $\ocG_K^\mathrm{rdt}$ is given by specifying linear automorphisms $A_{123}$ of  $\langle 1,2,3\rangle$ and $A_{45}$ of $\langle 4,5\rangle$ and a symplectic similitude $A_{\pm6,\pm7,\pm 8}$ of $\langle\pm 6,\pm 7,\pm 8\rangle$, such that $\left.A\right|_{\langle 1,2,3\rangle}=A_{123}$, $\left.A\right|_{\langle 4,5\rangle}=A_{45}$, $\left.A\right|_{\langle \pm 6,\pm 7,\pm8\rangle}=A_{\pm 6,\pm 7,\pm 8}$, where $\left.A\right|_{\langle -1,-2,-3\rangle}$ is uniquely determined by $A_{123}$, $c(A_{\pm 6,\pm 7,\pm 8})$ ($c$ being the multiplier character) and the imposition that $A$ be a symplectic similitude, and similarly for $\left.A\right|_{\langle -4,-5\rangle}$.

    If for example we consider $J/2g\Z=\{0,\pm 2, \pm 3, \pm 5\}$ instead, we expect something along the lines of $\GSp(6)\times \GL(2)\times \GL(2)$ and indeed we additionally get the subbundles
    \begin{align*}
      \{-2, -1\},
      \{1, 2\},
      \{-3\},
      \{-5, -4, -3\}, \\
      \{3, 4, 5\},
      \{3\},
      \{3, 4, 5, 6, 7, 8, -8, -7, -6, -5, -4, -3\},
      \{1, 2, -2, -1\}.
    \end{align*}

    So an element $A$ of $\ocG_K^\mathrm{rdt}$ is given by specifying linear automorphisms $A_{12}$ of  $\langle 1,2\rangle$ and $A_{45}$ of $\langle 4,5\rangle$ and a symplectic similitude $A_{\pm6,\pm7,\pm 8}$ of $\langle\pm 6,\pm 7,\pm 8\rangle$ in a similar way to above.
  \end{itemize}
\end{Example}

\subsubsection{\texorpdfstring{$\ocG_K\EKORZip$}{GK-EKORZip} in the Siegel case}
\label{sec:ekorzip-siegel}

Recall that we denote the unipotent radical of $\ocG_K$ by $R_u\ocG_K$.

We divide out of $\ocG_K\AdmZip^\sim$ the action of the smooth normal subgroup $R_u\ocG_K\times R_u\ocG_K\subseteq\ocG_K\times\ocG_K$ and observe that $\ocG_K\times\ocG_K$ still acts on $[R_u\ocG_K\times R_u\ocG_K \backslash \ocG_K\AdmZip^\sim]=:\ocG_K\EKORZip^\sim$ (not a scheme).

We also define $\ocG_K\EKORZip:=[(\Delta({\ocG_K})\cdot(R_u\ocG_K\times R_u\ocG_K))\backslash \ocG_K\AdmZip^\sim]$.

\begin{Proposition}\label{orbit-rdt}
  We have well-defined morphisms
  \begin{alignat*}{2}
    (\ocG_K\times\ocG_K)/\tE_w &\to (\ocG_K^\mathrm{rdt}\times\ocG_K^\mathrm{rdt})/E_{\mathcal{Z}_w}, &\quad(X,Y)&\mapsto (X^\mathrm{rdt},Y^\mathrm{rdt}), \\
    \ocG_K/\tE_w &\to \ocG_K^\mathrm{rdt}/E_{\mathcal{Z}_w}, &\quad X&\mapsto X^\mathrm{rdt},
  \end{alignat*}
  and a bijectiion
  \begin{equation*}
    (\ocG_K\times\ocG_K)/(\tE_w\cdot(R_u\ocG_K\times R_u\ocG_K)) \to (\ocG_K^\mathrm{rdt}\times\ocG_K^\mathrm{rdt})/E_{\mathcal{Z}_w}.
  \end{equation*}
\end{Proposition}

\begin{Proof}
  The first assertion follows from the definition of $E_{\mathcal{Z}_w}$ and equation~\eqref{eq:verallg-zip-gp}. The second then follows from Lemma~\ref{Ebild}.
\end{Proof}

\begin{Lemma}Assume $0\in J$.  
  The underlying topological spaces of the stacks in consideration are as follows:
  \begin{enumerate}[(1)]
  \item\label{item:stack-pts-1} $|[\ocG_K\backslash \oM^\mathrm{loc}]| = \KR(K,\{\mu\}) \overset{\text{def.}}= W_K\backslash (W_K \Adm(\{\mu\})W_K)/W_K$.
  \item\label{item:stack-pts-3} $|\ocG_K\EKORZip| = \EKOR(K,\{\mu\}) = \Adm(\{\mu\})^K\cap{}^K\tW$ \\
    $\overset{\text{\ref{iw-diagr}}}\cong \bigcup_{w\in \KR(K,\{\mu\})} \bK w \bK / \bK_\sigma(\bK_1\times \bK_1)$.
  \end{enumerate}
\end{Lemma}

\begin{Proof}
  \eqref{item:stack-pts-1} is well-known as explained in Section~\ref{sec:kottw-rapop-strat}.

  \eqref{item:stack-pts-3}: By Lemma~\ref{orbits-kr}, the $\ocG_K\times\ocG_K$-orbits in $\ocG_K\AdmZip^{\sim}$ are indexed by $\Adm(\{\mu\})_K=\KR(K,\{\mu\})$.

  Let us further investigate the $\ocG_K\times\ocG_K$-orbit of $\tStd(w)$ in $\ocG_K\EKORZip^{\sim}$ for some fixed $w\in\Adm(\{\mu\})^K$.

  By Proposition~\ref{orbit-rdt}, its underlying topological space agrees with that of $\ocG_K^\mathrm{rdt}\Zip^{\sim,\mathcal{Z}_w}$. By \cite{SYZnew} we know that $|\ocG_K^\mathrm{rdt}\Zip^{\mathcal{Z}_w}|\cong \bK w \bK / \bK_\sigma(\bK_1\times \bK_1)$, whence the lemma.
\end{Proof}

\begin{Corollary}
  We have a morphism
  \begin{equation*}
    \left(\ocG_K\AdmZip_w\right)_\mathrm{red}=\text{orbit of } \Std(w) \to \ocG_K^\mathrm{rdt}\Zip^{\mathcal{Z}_w}.
  \end{equation*}

  This defines the EKOR stratification on $\ocG_K\AdmZip_w$.  All in all, we get an EKOR stratification on $\ocG_K\AdmZip$.

  The morphism factors through $(\ocG_K\EKORZip_w)_\mathrm{red}$, and $(\ocG_K\EKORZip_w)_\mathrm{red}\to \ocG_K^\mathrm{rdt}\Zip^{\mathcal{Z}_w}$ is an isomorphism.
\end{Corollary}

\begin{Corollary}
  For every point of $[\ocG_K\backslash\oM_K]$, $\osS_K\to \ocG_K\EKORZip$ is smooth as a map between the associated reduced fiber of $\osS_K \to [\ocG_K\backslash\oM_K]$ and the associated reduced fiber of $\ocG_K\EKORZip \to [\ocG_K\backslash\oM_K]$.
\end{Corollary}

\begin{Proof}
  This follows from the preceding corollary by \cite[Theorem~A]{SYZnew} (which says that the map $\osS_{K}^w\to\ocG_K^\mathrm{rdt}\Zip^{\mathcal{Z}_w}$ is smooth, cf. subsection~\ref{sec:ekor-stratification}).
\end{Proof}

The key obstacle in going forward toward proving smoothness of $\osS_K \to \ocG_K\EKORZip$ now is that
we do not know whether the fibers of $\ocG_K\EKORZip \to [\ocG_K\backslash\oM_K]$ are reduced.

\begin{Conjecture}
  We conjecture that the answer is affirmative. In fact, we conjecture that $\ocG_K\EKORZip \to [\ocG_K\backslash\oM_K]$ is smooth.
\end{Conjecture}

\begin{Corollary}
  $\osS_K\to\ocG_K\EKORZip$ is surjective.
\end{Corollary}

\begin{Proof}
  This follows from the description of the topological space and what is already known from \cite[first paragraph of section~6.3]{he-rapo}.
\end{Proof}

We get a commutative diagram
\begin{equation*}
  \begin{tikzpicture}[node distance=4cm, auto]
    \node (LM) at (0,3) {$\tosS_K$};
    \node (RM) at (3,3) {$\ocG_K\AdmZip^\sim$};
    \node (RRM) at (7.1,3) {$\ocG_K\EKORZip^\sim$};
    \node (LU) at (0,0) {$\osS_K$};
    \node (RU) at (3,0) {$\ocG_K\AdmZip$};
    \node (RRU) at (7.1,0) {$\ocG_K\EKORZip$};
    \node (RRRU) at (11,0) {$[\ocG_K\backslash\oM^\mathrm{loc}_K]$};
    
    \draw[->] (LM) to  (LU);
    \draw[->] (RM) to  (RU);
    \draw[->] (RRM) to  (RRU);
    \draw[->] (LM) to  (RM);
    \draw[->] (LU) to  (RU);
    \draw[->] (RM) to (RRM);
    \draw[->] (RU) to  (RRU);
    \draw[->] (RRU) to  (RRRU);
  \end{tikzpicture}
\end{equation*}

\begin{Remark}\label{redRusmooth}
  Since $R_u\ocG_K$ is smooth, $\ocG_K\AdmZip^\sim \to \ocG_K\EKORZip^\sim$ is smooth.
\end{Remark}

\begin{Remark}
  Another open question at this point is: what is the relationship between $\ocG_K\text{-EKORZip}^\mathrm{perf}$ and the shtuka approach of \cite[Section~4]{SYZnew}?
\end{Remark}

\begin{Remark}
  It should be straightforward to generalize (taking into account the extra structure) our constructions to those (P)EL cases where the local model is the ``naive'' local model of Rapoport-Zink \cite{rz}.
\end{Remark}

\subsubsection{The example of \texorpdfstring{$\GSp(4)$}{GSp(4)}}
\label{sec:example-gsp4}

To illustrate some aspects, we look at the example $2g=4$. 

\paragraph{The apartment.}
\label{sec:apartment}

We describe the (extended) apartment. We follow the general outline of \cite{landvogt}, in particular as far as notation is concerned.

The roots are $\pm (2e_1-e_0), \pm(2e_2-e_0),\pm(e_1-e_2),\pm(e_1+e_2-e_0)$. The simple affine roots and the (various variants of the) Weyl group are as described in Remark~\ref{root-stuff-gsp}. The root one-parameter subgroups\footnote{The parameter being additive here; i.e., we're talking about homomorphisms $\G_a\to G$.} are given as follows:
\begin{align*}
  u_{e_1-e_2}(x) &=
                   \begin{pmatrix}
                     1 & x & &  \\
                     & 1 & &  \\
                     & & 1 & -x  \\
                     & & & 1
                   \end{pmatrix}, &
  u_{e_2-e_1}(x) &=
                   \begin{pmatrix}
                     1 &  & &  \\
                     x & 1 & &  \\
                     & & 1 &   \\
                     & & -x & 1
                   \end{pmatrix}, \displaybreak[0]\\
  u_{2e_1-e_0}(x) &=
                   \begin{pmatrix}
                     1 & & & x  \\
                     & 1 & &  \\
                     & & 1 &  \\
                     & & & 1
                   \end{pmatrix}, &
  u_{e_0-2e_1}(x) &=
                   \begin{pmatrix}
                     1 & & &   \\
                     & 1 & &  \\
                     & & 1 &  \\
                     x & & & 1
                   \end{pmatrix}, \displaybreak[0]\\
  u_{2e_2-e_0}(x) &=
                   \begin{pmatrix}
                     1 & & &  \\
                     & 1 & x &  \\
                     & & 1 &  \\
                     & & & 1
                   \end{pmatrix}, &
  u_{e_0-2e_2}(x) &=
                   \begin{pmatrix}
                     1 & & &  \\
                     & 1 & &  \\
                     & x & 1 &  \\
                     & & & 1
                   \end{pmatrix}, \displaybreak[0]\\
  u_{e_1+e_2-e_0}(x) &=
                   \begin{pmatrix}
                     1 &  &x &  \\
                     & 1 & & x  \\
                     & & 1 &  \\
                     & & & 1
                   \end{pmatrix}, &
  u_{e_0-e_1-e_2}(x) &=
                   \begin{pmatrix}
                     1 &  & &  \\
                     & 1 & &   \\
                     x & & 1 &  \\
                     & x & & 1
                   \end{pmatrix}
\end{align*}

For $a\in R$ define $w_a(x):=u_a(x)u_{-a}(-x^{-1})u_a(x)$.

\begin{Remark}\label{gsp4-normalizer}
  $N(\Q_p)$ is generated by $T(\Q_p)$ and all $w_a(x)$ as above.
\end{Remark}

\begin{Remark}
  $w_a(x)=m(u_{-a}(-x^{-1}))$ in Landvogt's notation \cite{landvogt}.
\end{Remark}

We have $V_1:=X_*(T)\otimes\R=\{(x_1,x_2,x_{-2},x_{-1})\in\R^4 \suchthat x_1+x_{-1}=x_2+x_{-2}\}$ and
\begin{equation*}
  \nu_1\colon T(\Q_p) \to V_1, \;
  \begin{pmatrix}
  d_1 & & & \\
  & d_2 & & \\
  & & cd_2^{-1} & \\
  & & & cd_1^{-1}
\end{pmatrix} \mapsto
\begin{pmatrix}
  -v_p(d_1) \\
  -v_p(d_2) \\
  -v_p(cd_2^{-1}) \\
  -v_p(cd_1^{-1})
\end{pmatrix}.
\end{equation*}

Also, $V_0=\{v\in V_1 \suchthat a(v)=0 \; \forall a\in \Phi\}=\R(1,1,1,1)$, $V:=V_1/V_0$.

The extended apartment $A=A^\mathrm{ext}$ now is an affine $V_1$-space together with the map $\nu_1\colon N(\Q_p)\to \Aff(A)=\GL(V_1)\ltimes V_1$, whose restriction to $T(\Q_p)$ is given as above and (cf.~Remark~\ref{gsp4-normalizer})
\begin{align*}
  \nu_1(w_{2e_1-e_0}(x)) &= (
  \begin{smallpmatrix}
    & & & 1 \\
    & 1 & & \\
    & & 1 & \\
    1 & & &
  \end{smallpmatrix},
  \begin{smallpmatrix}
    -v_p(x) \\
    0 \\
    0 \\
    v_p(x)
  \end{smallpmatrix}
  ), \displaybreak[0]\\
  \nu_1(w_{2e_2-e_0}(x)) &= (
  \begin{smallpmatrix}
    1 & & &  \\
    &  & 1 & \\
    & 1 &  & \\
     & & & 1
  \end{smallpmatrix},
  \begin{smallpmatrix}
    0 \\
    -v_p(x) \\
    v_p(x) \\
    0
  \end{smallpmatrix}
  ), \displaybreak[0]\\
  \nu_1(w_{e_1-e_2}(x)) &= (
  \begin{smallpmatrix}
    & 1 & &  \\
    1 &  & & \\
    & &  & 1\\
    & & 1 &
  \end{smallpmatrix},
  \begin{smallpmatrix}
    -v_p(x) \\
    v_p(x) \\
    -v_p(x) \\
    v_p(x)
  \end{smallpmatrix}
  ), \displaybreak[0]\\
  \nu_1(w_{e_1+e_2-e_0}(x)) &= (
  \begin{smallpmatrix}
    & & 1 &  \\
    &  & & 1 \\
    1 & &  & \\
     & 1 & &
  \end{smallpmatrix},
  \begin{smallpmatrix}
    -v_p(x) \\
    -v_p(x) \\
    v_p(x) \\
    v_p(x)
  \end{smallpmatrix}
  ),
\end{align*}
etc. (Recipe: Write $w_a(x)$ as a product of a diagonal matrix $\diag(d_1,d_2,d_{-2},d_{-1})$ and a permutation matrix $P$ (this need not be a factorization in $\GSp(4)$); then
\begin{equation*}
\nu_1(w_a(x))=(P,
\begin{smallpmatrix}
  -v_p(d_1) \\
  -v_p(d_2) \\
  -v_p(d_{-2}) \\
  -v_p(d_{-1})
\end{smallpmatrix}
).)
\end{equation*}

The reduced apartment $A^\mathrm{red}$ is the affine $V$-space together with $\nu\colon N(\Q_p)\to \Aff(A^\mathrm{red})=\GL(V)\ltimes V$ given by the same formulas.

The walls (or rather, wall conditions) are given as follows ($n\in\Z$):
\begin{align*}
  2e_1-e_0 &: n=x_0-2x_1, \\
  2e_2-e_0 &: n=x_0-2x_2, \\
  e_1-e_2 &: n=x_2-x_1, \\
  e_1+e_2-e_0 &: n=x_0-x_1-x_2.
\end{align*}

\begin{figure}
  \caption{The reduced apartment with the base alcove highlighted.}\label{fig:red-apt}
  \centering
  \includegraphics[width=.7\textwidth]{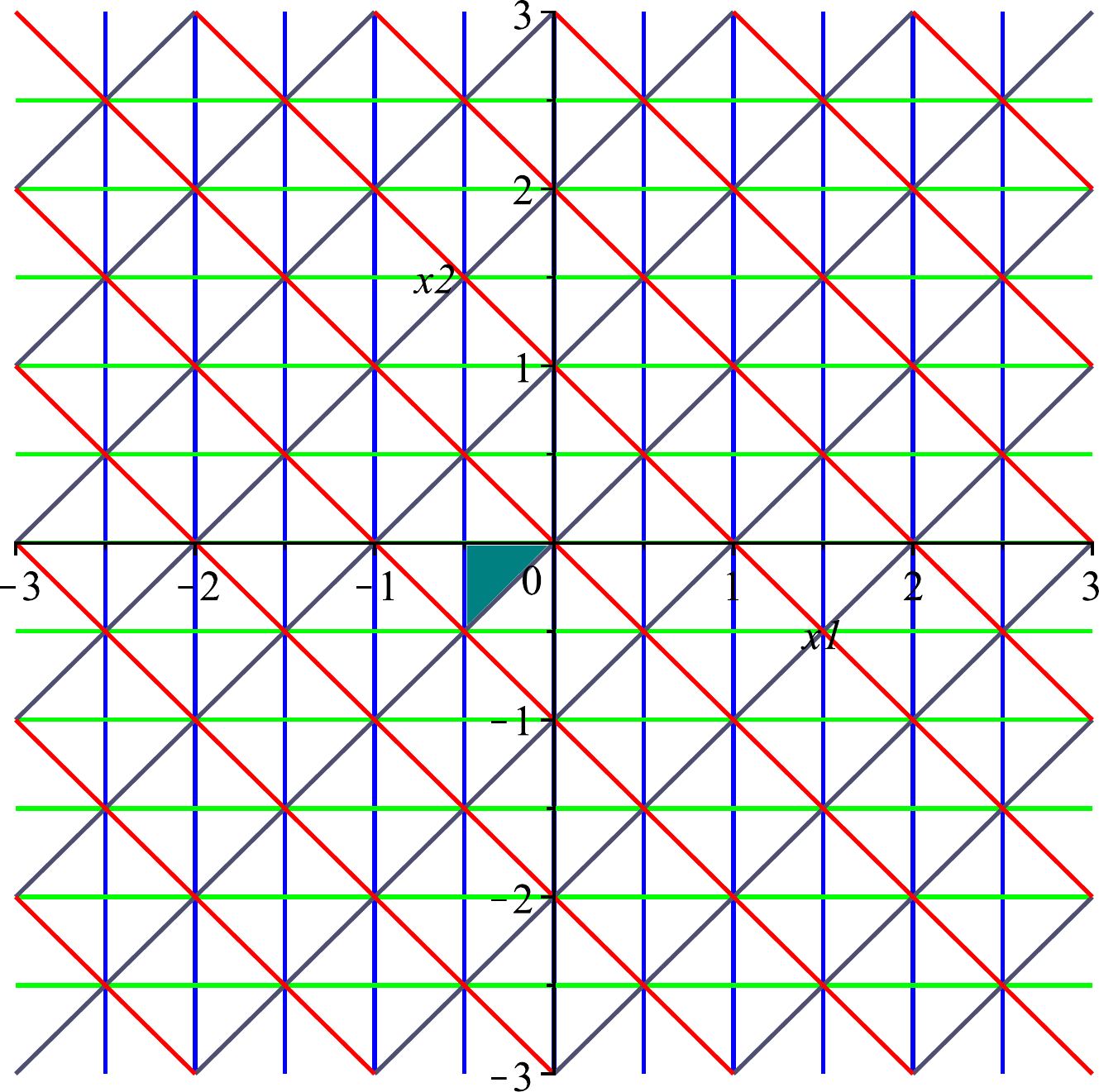}
\end{figure}

\paragraph{Lattice chains and parahoric subgroups.}
\label{sec:lattice-chains}

By \cite{zbMATH03900941}, the extended building $\mathcal{B}(\GL(X),\Q_p)$ is in bijection with norms\footnote{Defining conditions for a norm: $\alpha(tx)=\alpha(x)+\ord_p(t)$, $\alpha(x+y)\geq \min(\alpha(x),\alpha(y))$, $\alpha(x)=\infty \iff x=0$} $\alpha\colon X\to \R\cup \{\infty\}$. Norms in turn are in bijection with graded lattice chains (cf. Remark~\ref{gl-building}). Indeed, if $\alpha$ is a norm, define $\Delta_\alpha$ to be the set of its balls centered around zero and $c_\alpha(\Lambda):=\inf_{\lambda\in\Lambda}\alpha(\lambda)$. Conversely, given a graded lattice chain $(\Delta,c)$, define a norm $\alpha$ by $\alpha(x):=c(\Lambda)$ for the smallest $\Lambda\in\Delta$ with $x\in\Lambda$.

To go from the extended apartment of $\GL(X)$, an affine $\R^{n}$-space, where $n=\dim X$, to norms, fix a basis $e_1,\dotsc,e_n$ of $X$. Then $v\in\R^n$ corresponds to the norm $\alpha_v$ with
\begin{equation*}
  \alpha_v(\sum t_ie_i)=\min_i(\ord_p(t_i) - v_i).
\end{equation*}

There are seven types of points in the extended apartment (in each case we choose one in the base alcove to represent all of its type) corresponding to the vertices, edges and interior of the base alcove:
\begin{itemize}
\item standard hyperspecial: $x_\mathrm{hs}=(0,0,0,0)$
\item paramodular: $x_\mathrm{paramod}=(-1/2,0,0,1/2)$
\item Klingen: $x_\mathrm{Klingen}=(-1/4,0,0,1/4)$
\item Siegel: $x_\mathrm{Siegel}=(-1/4,-1/4,1/4,1/4)$
\item Iwahori: $x_\mathrm{Iwahori}=(-1/4,-1/8,1/8,1/4)$
\item another hyperspecial: $x=(-1/2,-1/2,1/2,1/2)$
\item another parahoric: $x=(-1/2,-1/4,1/4,1/2)$
\end{itemize}

The last two are conjugates (by the Atkin-Lehner element) of the standard hyperspecial and the Klingen parahoric, respectively (see e.g. \cite[151]{roesner}); therefore we will neglect them in the sequel.

For a set of lattices $S$ denote by $\langle S \rangle$ the closure under homotheties, i.e., $\langle S\rangle :=\{ p^n s \suchthat n\in \Z,\; s\in S\}$.

Then:
\begin{itemize}
\item $\Delta_\mathrm{hs}=\langle \Z_p^4\rangle$ \\ and $c_\mathrm{hs}(\Z_p^4)=0$.
\item $\Delta_\mathrm{paramod}=\langle \Z_p^3\oplus p\Z_p, \; \Z_p\oplus p\Z_p^3\rangle$ \\ and $c_\mathrm{paramod}(\Z_p^3\oplus p\Z_p)=-\frac{1}{2}$, $c_\mathrm{paramod}(\Z_p\oplus p\Z_p^3)=0$.
\item $\Delta_\mathrm{Klingen}=\langle \Z_p^4, \Z_p^3\oplus p\Z_p,\Z_p\oplus p\Z_p^3\rangle$
\\ and $c_\mathrm{Klingen}(\Z_p^4)=-1/4$, $c_\mathrm{Klingen}(\Z_p^3\oplus p\Z_p)=0$, $c_\mathrm{Klingen}(\Z_p\oplus p\Z_p^3)=1/4$.
\item $\Delta_\mathrm{Siegel}=\langle \Z_p^4, \Z_p^2\oplus p\Z_p^2\rangle$
\\ and $c_\mathrm{Siegel}(\Z_p^4)=-1/4$, $c_\mathrm{Siegel}(\Z_p^2\oplus p\Z_p^2)=1/4$.
\item $\Delta_\mathrm{Iwahori}=\langle \Z_p^4, \Z_p^3\oplus p\Z_p,\Z_p^2\oplus p\Z_p^2,\Z_p\oplus p\Z_p^3\rangle$
\\ and $c_\mathrm{Iwahori}(\Z_p^4)=-1/4$, $c_\mathrm{Iwahori}(\Z_p^3\oplus p\Z_p)=-1/8$, $c_\mathrm{Iwahori}(\Z_p^2\oplus p\Z_p^2)=1/8$, $c_\mathrm{Iwahori}(\Z_p\oplus p\Z_p^3)=1/4$.
\end{itemize}

The associated parahoric subgroups are
\begin{itemize}
\item hyperspecial: $\GSp_4(\Z_p)$
\item paramodular: $\GSp_4(\Q_p) \cap
  \begin{pmatrix}
    \Z_p & \Z_p & \Z_p & p^{-1}\Z_p \\
    p\Z_p & \Z_p & \Z_p & \Z_p \\
    p\Z_p & \Z_p & \Z_p & \Z_p \\
    p\Z_p & p\Z_p & p\Z_p & \Z_p
  \end{pmatrix}
$
\item Klingen: $\GSp_4(\Z_p) \cap \begin{pmatrix}
    \Z_p & \Z_p & \Z_p & \Z_p \\
    p\Z_p & \Z_p & \Z_p & \Z_p \\
    p\Z_p & \Z_p & \Z_p & \Z_p \\
    p\Z_p & p\Z_p & p\Z_p & \Z_p
  \end{pmatrix}$
\item Siegel: $\GSp_4(\Z_p) \cap \begin{pmatrix}
    \Z_p & \Z_p & \Z_p & \Z_p \\
    \Z_p & \Z_p & \Z_p & \Z_p \\
    p\Z_p & p\Z_p & \Z_p & \Z_p \\
    p\Z_p & p\Z_p & \Z_p & \Z_p
  \end{pmatrix}$
\item Iwahori: $\GSp_4(\Z_p) \cap \begin{pmatrix}
    \Z_p & \Z_p & \Z_p & \Z_p \\
    p\Z_p & \Z_p & \Z_p & \Z_p \\
    p\Z_p & p\Z_p & \Z_p & \Z_p \\
    p\Z_p & p\Z_p & p\Z_p & \Z_p
  \end{pmatrix}$
\end{itemize}

\begin{Remark} Dualizing with respect to the symplectic form, we have
  \begin{align*}
    (\Z_p^4)^\vee&=\Z_p^4, &
                   (\Z_p\oplus p\Z_p^3)^\vee &= p^{-1}(\Z_p^3\oplus p\Z_p), \\
    (\Z_p^3\oplus p\Z_p)^\vee&=p^{-1}(\Z_p\oplus p\Z_p^3), &
                               (\Z_p^2\oplus p\Z_p^2)^\vee &= p^{-1}(\Z_p^2\oplus p\Z_p^2).
\end{align*}

\end{Remark}

\paragraph{Admissible set.}
\label{sec:admissible-set}

We compute the admissible set in the way outlined in Remark~\ref{kottrap}. The cocharacter $\mu$ is $(1,1,0,0)$.

We obtain
\begin{align*}
  \Adm(\{\mu\})=\bigl\{
  &\Bigl(\id                               , \begin{smallpmatrix}1 \\ 1 \\ 0 \\ 0\end{smallpmatrix}\Bigr), \quad
  \Bigl((2 \quad {-2})                    , \begin{smallpmatrix}1 \\ 0 \\ 1 \\ 0\end{smallpmatrix}\Bigr), \quad
\Bigl(\id                               , \begin{smallpmatrix}1 \\ 0 \\ 1 \\ 0\end{smallpmatrix}\Bigr), \\
  &\Bigl((1 \quad {-1})                    , \begin{smallpmatrix}0 \\ 1 \\ 0 \\ 1\end{smallpmatrix}\Bigr), \quad
  \Bigl((1 \quad 2 \quad {-1} \quad {-2}) , \begin{smallpmatrix}0 \\ 1 \\ 0 \\ 1\end{smallpmatrix}\Bigr), \\
&\Bigl((1 \quad 2)({-2} \quad {-1})      , \begin{smallpmatrix}0 \\ 1 \\ 0 \\ 1\end{smallpmatrix}\Bigr), \quad
  \Bigl(\id                               , \begin{smallpmatrix}0 \\ 1 \\ 0 \\ 1\end{smallpmatrix}\Bigr), \\
&\Bigl((1 \quad {-1})(2 \quad {-2})      , \begin{smallpmatrix}0 \\ 0 \\ 1 \\ 1\end{smallpmatrix}\Bigr), \quad
\Bigl((1 \quad {-1})                    , \begin{smallpmatrix}0 \\ 0 \\ 1 \\ 1\end{smallpmatrix}\Bigr), \\
  &\Bigl((1 \quad {-2})(2 \quad {-1})      , \begin{smallpmatrix}0 \\ 0 \\ 1 \\ 1\end{smallpmatrix}\Bigr), \quad
\Bigl((1 \quad {-2} \quad {-1} \quad 2) , \begin{smallpmatrix}0 \\ 0 \\ 1 \\ 1\end{smallpmatrix}\Bigr), \\
&\Bigl((2 \quad {-2})                    , \begin{smallpmatrix}0 \\ 0 \\ 1 \\ 1\end{smallpmatrix}\Bigr), \quad
  \Bigl(\id                               , \begin{smallpmatrix}0 \\ 0 \\ 1 \\ 1\end{smallpmatrix}\Bigr)
  \bigr\},
\end{align*}
or, in terms of Frobenii (cf. Construction~\ref{std-zeug})
\begin{align*}
  \Bigl\{
  &\begin{smallpmatrix}
p & 0 & 0 & 0 \\
0 & p & 0 & 0 \\
0 & 0 & 1 & 0 \\
0 & 0 & 0 & 1
\end{smallpmatrix},
\begin{smallpmatrix}
p & 0 & 0 & 0 \\
0 & 0 & 1 & 0 \\
0 & p & 0 & 0 \\
0 & 0 & 0 & 1
\end{smallpmatrix},
\begin{smallpmatrix}
p & 0 & 0 & 0 \\
0 & 1 & 0 & 0 \\
0 & 0 & p & 0 \\
0 & 0 & 0 & 1
\end{smallpmatrix},
\begin{smallpmatrix}
0 & 0 & 0 & 1 \\
0 & p & 0 & 0 \\
0 & 0 & 1 & 0 \\
p & 0 & 0 & 0
\end{smallpmatrix},
\begin{smallpmatrix}
0 & 0 & 1 & 0 \\
p & 0 & 0 & 0 \\
0 & 0 & 0 & 1 \\
0 & p & 0 & 0
\end{smallpmatrix}, \\
&\begin{smallpmatrix}
0 & 1 & 0 & 0 \\
p & 0 & 0 & 0 \\
0 & 0 & 0 & 1 \\
0 & 0 & p & 0
\end{smallpmatrix},
\begin{smallpmatrix}
1 & 0 & 0 & 0 \\
0 & p & 0 & 0 \\
0 & 0 & 1 & 0 \\
0 & 0 & 0 & p
\end{smallpmatrix},
\begin{smallpmatrix}
0 & 0 & 0 & 1 \\
0 & 0 & 1 & 0 \\
0 & p & 0 & 0 \\
p & 0 & 0 & 0
\end{smallpmatrix},
\begin{smallpmatrix}
0 & 0 & 0 & 1 \\
0 & 1 & 0 & 0 \\
0 & 0 & p & 0 \\
p & 0 & 0 & 0
\end{smallpmatrix},
\begin{smallpmatrix}
0 & 0 & 1 & 0 \\
0 & 0 & 0 & 1 \\
p & 0 & 0 & 0 \\
0 & p & 0 & 0
\end{smallpmatrix}, \\
&\begin{smallpmatrix}
0 & 1 & 0 & 0 \\
0 & 0 & 0 & 1 \\
p & 0 & 0 & 0 \\
0 & 0 & p & 0
\end{smallpmatrix},
\begin{smallpmatrix}
1 & 0 & 0 & 0 \\
0 & 0 & 1 & 0 \\
0 & p & 0 & 0 \\
0 & 0 & 0 & p
\end{smallpmatrix},
\begin{smallpmatrix}
1 & 0 & 0 & 0 \\
0 & 1 & 0 & 0 \\
0 & 0 & p & 0 \\
0 & 0 & 0 & p
\end{smallpmatrix}
            \Bigr\}.
\end{align*}

\paragraph{Siegel level.}
\label{sec:siegel-level}

From now on, we consider the Siegel level structure. Denote the Siegel parahoric by $K$ and the standard hyperspecial subgroup by $H$. Here $W_K$ is generated by $({-1}\quad{-2})(1\quad 2)$, while $W_H$ is generated by $W_K$ and $(2\quad{-2})$.

Recalling Remark~\ref{symplecticity}\,\eqref{item:symplecticity-2}, we note that one has a natural morphism $\ocG_K\Zip \to \ocG_H\Zip \times \ocG_H\Zip$.

We have
\begin{align*}
  \KR(K,\{\mu\}) &= \Bigl\{ %
                    (\id, \begin{smallpmatrix}1\\ 1\\ 0\\ 0\end{smallpmatrix}),
((1\quad {-2})(2\quad {-1}), \begin{smallpmatrix}0\\ 0\\ 1\\ 1\end{smallpmatrix}),
((1\quad 2\quad {-1}\quad {-2}), \begin{smallpmatrix}0\\ 1\\ 0\\ 1\end{smallpmatrix}), \\
&(\id, \begin{smallpmatrix}0\\ 0\\ 1\\ 1\end{smallpmatrix}),
((1\quad {-2} \quad {-1} \quad 2), \begin{smallpmatrix}0\\ 0\\ 1\\ 1\end{smallpmatrix}),
  ((1 \quad 2)({-2} \quad {-1}), \begin{smallpmatrix}0\\ 1\\ 0\\ 1\end{smallpmatrix}) 
                            \Bigr\}, \\
  \EKOR(K,\{\mu\}) &= \KR(K,\{\mu\}) \cup   %
                      \Bigl\{
                      (\id, \begin{smallpmatrix} 0 \\ 1 \\ 0 \\ 1\end{smallpmatrix}), \quad
                      ((2\quad{-2}), \begin{smallpmatrix}0 \\ 0 \\ 1 \\ 1 \end{smallpmatrix}), \quad
  ((1\quad{-1}), \begin{smallpmatrix}0 \\ 1 \\ 0 \\ 1 \end{smallpmatrix})
\Bigr\}.
\end{align*}

\newgeometry{margin=0.76cm} %
\begin{landscape}
  ~\vskip 1.5cm
  In the following table, $w^j$ is the isomorphism type of the $\ocG_H$-zip at position $j$. For $\mathcal{C}^\bullet, \mathcal{D}^\bullet$ we give (indices of) basis vectors. ``$\leftarrow$'' means ``same as in the column adjacent to the left''. $\alpha_0\colon \bar\F_p^4\to \bar\F_p^4$ is the projection onto the plane spanned by the $1,2$-coordinates, $\alpha_2$ the projection onto the plane spanned by the $-2,-1$-coordinates. By $\alpha_{j,\mathcal{C}^\bullet/\mathcal{D}^\bullet}$ we denote the induced maps on $\mathcal{V}^\bullet/\mathcal{C}^\bullet\oplus \mathcal{C}^\bullet$ and $\mathcal{D}^\bullet\oplus \mathcal{V}^\bullet/\mathcal{D}^\bullet$, respectively. Each $\mathcal{C}^j\subseteq \mathcal{V}^j$ has a canonical complement in terms of standard basis vectors. Importantly however, we will not always have a complementary \emph{chain} of linear subspaces. In any event, below we say what the $\alpha_{j,\mathcal{C}^\bullet/\mathcal{D}^\bullet}$ are the projection onto if interpreted as described. For instance, the projection onto $\emptyset$ is the zero map. So in that case $\mathcal{V}^\bullet/\mathcal{C}^\bullet\oplus \mathcal{C}^\bullet$ (or $\mathcal{D}^\bullet\oplus \mathcal{V}^\bullet/\mathcal{D}^\bullet$) is a chain of vector spaces with zero transition maps.

  \begin{table}[hbtp]
    \centering
\begin{tabular}{l||l|l|l|l|l|l|l |l|l|l|l}
  $w$ & KR-type & $\mathcal{C}^0$ & $\mathcal{D}^0$ & $\mathcal{C}^2$ & $\mathcal{D}^2$ & $w^0$ & $w^2$ & $\alpha_{2,\mathcal{C}^\bullet}$ & $\alpha_{0,\mathcal{C}^\bullet}$ & $\alpha_{2,\mathcal{D}^\bullet}$ & $\alpha_{0,\mathcal{D}^\bullet}$ \\
  \hline
  $(\id, \begin{smallpmatrix}1\\ 1\\ 0\\ 0\end{smallpmatrix})$ & $\leftarrow$ & $\{1, 2\}$ & $\{1, 2\}$ & $\{-2, -1\}$ & $\{-2, -1\}$ & $(-2\quad 1)({-1}\quad 2)$ & $\leftarrow$ & $\{1, 2\}$ & $\{-2, -1\}$ & $\{1, 2\}$ & $\{-2, -1\}$ \\
  $((1\quad {-2})(2\quad {-1}), \begin{smallpmatrix}0\\ 0\\ 1\\ 1\end{smallpmatrix})$ & $\leftarrow$ & $\{1, 2\}$ & $\{-2, -1\}$ & $\{1, 2\}$ & $\{-2, -1\}$ & $\id$ & $\leftarrow$ & $\emptyset$ & $\emptyset$ & $\emptyset$ & $\emptyset$ \\
  $((1\quad 2\quad {-1}\quad {-2}), \begin{smallpmatrix}0\\ 1\\ 0\\ 1\end{smallpmatrix})$ & $\leftarrow$ & $\{1, 2\}$ & $\{-2, 1\}$ & $\{-2, 1\}$ & $\{-2, -1\}$ & $(-2\quad 2)$ & $\leftarrow$ & $\{1\}$ & $\{-1\}$ & $\{2\}$ & $\{-2\}$ \\
  $(\id, \begin{smallpmatrix}0\\ 0\\ 1\\ 1\end{smallpmatrix})$ & $\leftarrow$ & $\{-2, -1\}$ & $\{-2, -1\}$ & $\{1, 2\}$ & $\{1, 2\}$ & $(-2\quad 1)({-1}\quad 2)$ & $\leftarrow$ & $\{1, 2\}$ & $\{-2, -1\}$ & $\{1, 2\}$ & $\{-2, -1\}$ \\
  $((1\quad {-2} \quad {-1} \quad 2), \begin{smallpmatrix}0\\ 0\\ 1\\ 1\end{smallpmatrix})$ & $\leftarrow$ & $\{-2, 1\}$ & $\{-2, -1\}$ & $\{1, 2\}$ & $\{-2, 1\}$ & $(-2\quad 2)$ & $\leftarrow$ & $\{2\}$ & $\{-2\}$ & $\{1\}$ & $\{-1\}$ \\
  $((1 \quad 2)({-2} \quad {-1}), \begin{smallpmatrix}0\\ 1\\ 0\\ 1\end{smallpmatrix})$ & $\leftarrow$ & $\{-2, 1\}$ & $\{-2, 1\}$ & $\{-2, 1\}$ & $\{-2, 1\}$ & $\id$ & $\leftarrow$ & $\{1, 2\}$ & $\{-2, -1\}$ & $\{1, 2\}$ & $\{-2, -1\}$ \\
  $(\id, \begin{smallpmatrix} 0 \\ 1 \\ 0 \\ 1\end{smallpmatrix})$ & $((1 \quad 2)({-2} \quad {-1}), \begin{smallpmatrix}0\\ 1\\ 0\\ 1\end{smallpmatrix})$ & $\{-1, 2\}$ & $\{-1, 2\}$ & $\{-2, 1\}$ & $\{-2, 1\}$ & $(-2\quad 1)({-1}\quad 2)$ & $\leftarrow$ & $\{1, 2\}$ & $\{-2, -1\}$ & $\{1, 2\}$ & $\{-2, -1\}$ \\
  $((2\quad{-2}), \begin{smallpmatrix}0 \\ 0 \\ 1 \\ 1 \end{smallpmatrix})$ & $((1\quad {-2} \quad {-1} \quad 2), \begin{smallpmatrix}0\\ 0\\ 1\\ 1\end{smallpmatrix})$ & $\{-1, 2\}$ & $\{-2, -1\}$ & $\{1, 2\}$ & $\{-2, 1\}$ & $(-2\quad{-1}\quad 2\quad 1)$ & $\leftarrow$ & $\{1\}$ & $\{-1\}$ & $\{1\}$ & $\{-1\}$ \\
  $((1\quad{-1}), \begin{smallpmatrix}0 \\ 1 \\ 0 \\ 1 \end{smallpmatrix})$ & $((1\quad 2\quad {-1}\quad {-2}), \begin{smallpmatrix}0\\ 1\\ 0\\ 1\end{smallpmatrix})$ & $\{1, 2\}$ & $\{-1, 2\}$ & $\{-2, 1\}$ & $\{-2, -1\}$ & $(-2\quad{-1}\quad 2\quad 1)$ & $\leftarrow$ & $\{2\}$ & $\{-2\}$ & $\{2\}$ & $\{-2\}$
\end{tabular}
\end{table}
\end{landscape}
\restoregeometry %

\begin{Observations}
  \begin{itemize}
  \item We always have $w^0=w^2$.  This is explained by the fact that the Ekedahl-Oort stratification in this case agrees with the Newton stratification (and isogenous abelian varieties by definition lie in the same Newton stratum).  
  \item Consider the Kottwitz-Rapoport strata containing more than one EKOR stratum (i.e., containing two EKOR strata). Then we can distinguish among the EKOR strata by looking at the Ekedahl-Oort stratum.  In other words, the EKOR stratification is in this case the coarsest common refinement of the Kottwitz-Rapoport and Ekedahl-Oort stratifications.
  \end{itemize}
\end{Observations}

\nocite{kottwitz,illusie,travaux,haines,zhangEO,zinka,zinkformal,lautrunc,rapo-guide,landvogt-buch,wortmann}

\printbibliography[heading=bibintoc]

\end{document}